\documentclass[12pt]{article}

\usepackage{myarticle}

\usepackage{mycommands}
\newcommand{\setsep}{\,\vert\,}
\newcommand{\Biggsetsep}{\,\Bigg\vert\,}

\newcommand{\ak}[1]{\bar{#1}}
\newcommand{\akk}[1]{\tilde{#1}}
\newcommand{\akp}[1]{\hat{#1}}

\newcommand{\X}{\R^N}                    % some space
\newcommand{\leg}{h}                      % legendre function
\newcommand{\clLeg}{\mathscr L}           % class of Legendre functions
\newcommand{\bregmap}[1][]{D_{#1}}        % Bregman distance (as a mapping without arguments)
\newcommand{\breg}[3][]{D_{#1}(#2,#3)}    % Bregman distance
\newcommand{\slp}[1]{\abs{\nabla #1}}     % slope of a function
\newcommand{\sslp}[1]{\overline{\abs{\nabla #1}}}    % limiting slope of a function
\newcommand{\pos}[1]{[#1]_+}              % max(#1,0) positive portion of a number

\newcommand{\fun}{f}                     % objective function
\newcommand{\mfun}[1]{f_{#1}}            % model function
\newcommand{\mfunp}[2]{f^{#1}_{#2}}      % proximized model function (f^#1_{#2} = f_#21 + 1/2\tau|x-#22|^2)
\newcommand{\mDeltap}[2]{\Delta^{#1}_{#2}}

\newcommand{\infbd}{\underline{\fun}}     % optimal value of the optimization problem

\newcommand{\kp}{{k+1}}
\renewcommand{\k}{{k}}

\newcommand{\iter}[1]{_{#1}}
\newcommand{\idx}[1]{^{(#1)}}

\newcommand{\Prox}[1][]{\ensuremath{P^{#1}}}            % proximal mapping

\newcommand{\landau}{\o}
\newcommand{\dir}{v}%\bm d}

\renewcommand{\O}{\mathcal O}

\colorlet{changecolor}{orange}

\usepackage{pgfplots}
\pgfplotsset{
  tick label style={font=\footnotesize},
  label style={font=\footnotesize},
  legend style={font=\footnotesize}
}

\graphicspath{{./figures/}}
\title{Non-smooth Non-convex Bregman Minimization:
Unification and New Algorithms}

\author{Peter Ochs, Thomas Brox, and Jalal Fadili
}


\begin{document}

\thispagestyle{empty}
\begin{center}
\vspace*{0.03\paperheight}
{\Large\bf\color{black}Non-smooth Non-convex Bregman Minimization:\\ \medskip
Unification and new Algorithms}\\
\bigskip
\bigskip
{\large Peter Ochs$^*$, Jalal Fadili$^\dagger$, and Thomas Brox$^\ddagger$ \\ \medskip
%\today 
{\small
$^*$~Saarland University, Saarbr\"{u}cken, Germany \\
$^\dagger$~Normandie Univ, ENSICAEN, CNRS, GREYC, France \\
$^\ddagger$~University of Freiburg, Freiburg, Germany \\
}
}
\end{center}
\bigskip

% ********************
% >>>>> ABSTRACT <<<<<
% ********************
\begin{center}
{\bf\color{black}Abstract}\\ \bigskip
\begin{minipage}{0.9\linewidth} \small
  We propose a unifying algorithm for non-smooth non-convex optimization. The algorithm approximates the objective function by a convex model function and finds an approximate (Bregman) proximal point of the convex model. This approximate minimizer of the model function yields a descent direction, along which the next iterate is found. Complemented with an Armijo-like line search strategy, we obtain a flexible algorithm for which we prove (subsequential) convergence to a stationary point under weak assumptions on the growth of the model function error. Special instances of the algorithm with a Euclidean distance function are, for example, Gradient Descent, Forward--Backward Splitting, ProxDescent, without the common requirement of a \enquote{Lipschitz continuous gradient}. In addition, we consider a broad class of Bregman distance functions (generated by Legendre functions), replacing the Euclidean distance.  The algorithm has a wide range of applications including many linear and non-linear inverse problems in signal/image processing and machine learning.
\end{minipage}
\end{center}
\bigskip

%\makekeywords

%\tableofcontents

% *******************
% >>>>> SECTION <<<<<
% *******************
\section{Introduction}

When minimizing a non-linear function on the Euclidean vector space, a fundamental strategy is to successively minimize approximations to the actual objective function. We refer to such an approximation as model (function). A common model example in smooth optimization is linearization (first order Taylor approximation) around the current iterate. However, in general, the minimization of a linear function does not provide a finite solution, unless, for instance, the domain is compact. Therefore, the model is usually complemented by a proximity measure, which favors a solution (the next iterate) close to the current iterate. For the Euclidean norm as proximity measure, computing the next iterate (minimizer of the sum of the model function and the Euclidean proximity measure) is equivalent to a Gradient Descent step, i.e. the next iterate is obtained by moving along the direction of the negative gradient at the current point for a certain step size. 

Since sequential minimization of model functions does not require smoothness of the objective or the model function, non-smoothness is handled naturally. The crucial aspect is the \enquote{approximation quality} of the model function, which is controlled by a growth function, that describes the approximation error around the current iterate. Drusvyatskiy et al. \cite{DIL16} refer to such model functions as Taylor-like models. The difference among algorithms lies in the properties of such a growth function, rather than the specific choice of a model function. \\

For the example of the Gradient Descent model function (linearization around the current iterate) for a continuously differentiable function, the value and the derivative of the growth function (approximation error) vanish at the current iterate. In this case, a line search strategy is required to determine a suitable step size that reduces the objective value. If the gradient of the objective function is additionally $L$-Lipschitz continuous, then the growth function satisfies a quadratic growth globally, and step sizes can be controlled analytically.\\

A large class of algorithms, which are widely popular in machine learning, statistics, computer vision, signal and image processing\xspace can be cast in the same framework. This includes algorithms such as Forward--Backward Splitting \cite{LM79} (Proximal Gradient Descent), ProxDescent \cite{LW16,DL16} (or proximal Gauss--Newton method), and many others. They all obey the same growth function as Gradient Descent. This allows for a unified analysis of all these algorithms, which is a key contribution of this paper. Moreover, we allow for a broad class of (iteration dependent) Bregman proximity functions (e.g., generated by common entropies such as Boltzmann--Shannon, Fermi--Dirac, and Burg's entropy), which leads to new algorithms. To be generic in the choice of the objective, the model, and the Bregman functions, the algorithm is complemented with an Armijo-like line search strategy. Subsequential convergence to a stationary point is established for different types of growth functions.\\

The above mentioned algorithms are ubiquitous in applications of machine learning, computer vision, image/signal processing, and statistics as is illustrated in Section~\ref{sec:examples} and in our numerical experiments in Section~\ref{sec:numerics}. Due to the unifying framework, the flexibility of these methods is considerably increased further.

% *******************
% >>>>> SECTION <<<<<
% *******************
\section{Contributions and Related Work} 

For smooth functions, Taylor's approximation is unique. However, for non-smooth functions, there are only \enquote{Taylor-like} model functions \cite{NPA08,Noll13,DIL16}. Each model function yields another algorithm. Some model functions \cite{NPA08,Noll13} could also be referred to as lower-Taylor-like models, as there is only a lower bound on the approximation quality of the model. Noll et al. \cite{Noll13} addressed the problem by bundle methods based on cutting planes, which differs from our setup.\\

The goal of Drusvyatskiy et al. \cite{DIL16} is to measure the proximity of an approximate solution of the model function to a stationary point of the original objective, i.e., a suitable stopping criterion for non-smooth objectives is sought. On the one hand, their model functions may be non-convex, unlike ours. On the other hand, their growth functions are more restrictive. Considering their abstract level, the convergence results may seem satisfactory. However, several assumptions that do not allow for a concrete implementation are required, such as a vanishing distance between successive iterates and convergence of the objective values along a generated convergent subsequence to the objective value of the limit point. This is in contrast to our framework. \\

We assume more structure of the subproblems: They are given as the sum of a model function and a Bregman proximity function. With this mild assumption on the structure and a suitable line-search procedure, the algorithm can be implemented and the convergence results apply. \emph{We present the first implementable algorithm in the abstract model function framework} and \emph{prove subsequential convergence to a stationary point}.\\

Our algorithm \emph{generalizes ProxDescent \cite{DL16,LW16} with convex subproblems}, which is known for its broad applicability. We provide more flexibility by considering Bregman proximity functions, and our \emph{backtracking line-search need not solve the subproblems for each trial step}.\\

The algorithm and convergence analysis is a \emph{far-reaching generalization} of Bonettini et al. \cite{BLPP16}, which is similar to the instantiation of our framework where the model function leads to Forward--Backward Splitting. The proximity measure of Bonettini et al. \cite{BLPP16} is assumed to satisfy a strong convexity assumption. Our \emph{proximity functions can be generated by a broad class of Legendre functions}, which includes, for example, the non-strongly convex Burg's entropy \cite{Burg72,BB97} for the generation of the Bregman proximity function.

% *******************
% >>>>> SECTION <<<<<
% *******************
\section{Preliminaries and Notations} \label{sec:prelim}

Throughout the whole paper, we work in a Euclidean vector space $\X$ of dimension $N\in\N$ equipped with the standard \emphdef{inner product} $\scal\cdot\cdot$ and associated \emphdef{norm} $\vnorm\cdot$. 

\paragraph{Variational analysis.} We work with extended-valued functions $\map{f}{\X}{\eR}$, $\eR:= \R\cup\set{\pm\infty}$. The \emphdef{domain} of $f$ is $\dom f:= \set{x\in \X\setsep f(x) < +\infty}$ and a function $f$ is \emphdef{proper}, if it is nowhere $-\infty$ and $\dom f\neq\emptyset$. It is \emphdef{lower semi-continuous} (or \emphdef{closed}), if $\liminf_{x\to\bar x} f(x) \geq f(\bar x)$ for any $\bar x\in \X$. Let $\sint \Omega$ denote the \emphdef{interior} of $\Omega\subset\X$. We use the notation of \emphdef{$f$-attentive convergence} $x \fto f \bar x \Leftrightarrow (x,f(x)) \to (\bar x, f(\bar x))$, and the notation $\k\fto K\infty$ for some $K\subset\N$ to represent $\k\to\infty$ where $\k\in K$.

As in \cite{DIL16}, we introduce the following concepts. For a closed function $\map{f}{\X}{\eR}$ and a point $\bar x \in \dom f$, we define the \emphdef{slope} of $f$ at $\bar x$ by 
\[
  \slp f (\bar x) := \limsup_{{x\to \bar x,\, x\neq \bar x}}\, \frac{\pos{f(\bar x) - f(x)}}{\vnorm{x-\bar x}}\,,
\]
where $\pos s := \max\{s,0\}$. It is the maximal instantaneous rate of decrease of $f$ at $\bar x$. For a differentiable function, it coincides with the norm of the gradient $\vnorm{\nabla f(\bar x)}$. Moreover, the \emphdef{limiting slope} 
  \[
    \sslp f (\bar x) := \liminf_{x\fto f \bar x}\, \slp f(x)
  \]
  is key. For a convex function $f$, we have $\sslp f (\bar x)= \inf_{v\in\partial f(\bar x)}\abs{v}$, where $\partial f(\bar x)$ is the \emphdef{(convex) subdifferential} $\partial f(\bar x) :=\set{v\in\R^N\setsep \forall x\colon f(x) \geq f(\bar x) + \scal{x-\bar x}{v}}$, whose \emphdef{domain} is given by $\dom\partial f:=\set{x\in\R^N\setsep \partial f(x) \neq \emptyset}$.  A point $\bar x$ is a \emphdef{stationary point} of the function $f$, if $\sslp f(\bar x) = 0$ holds. Obviously, if $\slp f(\bar x) = 0$, then $\sslp f(\bar x) = 0$. We define the \emphdef{set of (global) minimizers} of a function $f$ by
\[
    \Argmin_{x\in \X}\, f(x) := \set{x\in \X\setsep f(x) = \inf_{\bar x\in\X} f(\bar x) }\,,
\]
and the \emphdef{(unique) minimizer} of $f$ by $\argmin_{x\in\X}\, f(x)$, if $\Argmin_{x\in \X}\, f(x)$ consists of a single element. As shorthand, we also use $\Argmin f$ and $\argmin f$.

\begin{DEF}[{Growth function \cite{DIL16}}] \label{def:growth-function} 
  A differentiable univariate function $\map{\omega}{\R_+}{\R_+}$ is called \emphdef{growth function} if it satisfies $\omega(0)=\omega_+^\prime(0) = 0$. If, in addition, $\omega_+^\prime(t) >0$ for $t>0$ and\xspace equalities $\lim_{t\dto 0} \omega_+^\prime(t) = \lim_{t\dto 0} \omega(t)/\omega_+^\prime(t) = 0$ hold, we say that $\omega$ is a \emphdef{proper growth function}.
\end{DEF}
Concrete instances of growth functions will be generated for example by the concept of $\psi$-uniform continuity, which is a generalization of Lipschitz and H\"older continuity.
\begin{DEF}
  A mapping $\map{F}{\R^N}{\R^M}$ is called $\psi$-uniform continuous with respect to a continuous function $\map{\psi}{\R_+}{\R_+}$ with $\psi(0)=0$, if the following holds:
  \[
    \vnorm{F(x) - F(\bar x)} \leq \psi(\vnorm{x- \bar x})\quad \text{for all}\ x,\bar x\in \R^N \,.
  \]
\end{DEF}
\begin{EX}
  Let $F$ be $\psi$-uniform continuous. If, for some $c>0$, we have $\psi(s)=cs^\alpha$ with $\alpha\in ]0,1]$, then $F$ is H\"older continuous, which for $\alpha=1$ is the same as Lipschitz continuity.
\end{EX}
In analogy to the case of Lipschitz continuity, we can state a generalized Descent Lemma:
\begin{LEM}[Generalized Descent Lemma] \label{lem:gen-descent-lemma}
  Let $\map{f}{\R^N}{\R}$ be continuously differentiable and let $\map{\nabla f}{\R^N}{\R^N}$ be $\psi$-uniform continuous. Then, the following holds
  \[
    \abs{ f(x) - f(\bar x) - \scal{\nabla f(\bar x)}{x-\bar x} } \leq \int_0^1 \frac{\phi(s\vnorm{x-\bar x})}{s}\, \mathit{d}s \quad \text{for all}\  x,\bar x\in \R^N\,,
  \]
  where $\map{\phi}{\R_+}{\R_+}$ is given by $\phi(s):=s\psi(s)$.
\end{LEM}
\begin{proof}
  We follow the proof of the Descent Lemma for functions with Lipschitz gradient:
  \[
    \begin{split}
    \abs{f(x) - f(\bar x) - \scal{\nabla f(\bar x)}{x-\bar x}} 
    =&\ \abs{ \int_0^1 {\scal{\nabla f(\bar x+s(x-\bar x))-\nabla f(\bar x)}{x-\bar x}}\, \mathit{d}s } \\
    \leq&\  \int_0^1 {\vnorm{\nabla f(\bar x+s(x-\bar x))-\nabla f(\bar x)}\vnorm{x-\bar x}}\, \mathit{d}s \\
    \leq&\ \int_0^1 {\psi(s\vnorm{x-\bar x})\vnorm{x-\bar x}}\, \mathit{d}s  
    = \int_0^1 \frac{\phi(s\vnorm{x-\bar x})}{s}\, \mathit{d}s \,.
    \qedhere
    \end{split}
  \]
\end{proof}
\begin{EX}
The function $\omega(t)=\int_0^1 \frac{\phi(st)}{s}\, \mathit{d}s$ is an example for a growth function. Obviously, we have $\omega(0)=0$ and, using the Dominated Convergence Theorem (with majorizer $\sup_{s\in[0,1]}\psi(s)<+\infty$ for small $t\geq 0$), we conclude 
\[
  \omega_+^\prime(0)
    =\lim_{t\dto0} \int_0^1 \frac{\phi(st)}{st}\, \mathit{d}s   
    =\lim_{t\dto0} \int_0^1 {\psi(st)}\, \mathit{d}s   
    = \int_0^1 {\lim_{t\dto0} \psi(st)} \, \mathit{d}s = 0 \,.
\]
It becomes a proper growth function, for example, if $\psi(s)=0\Leftrightarrow s=0$ and we impose the additional condition $\lim_{t \dto 0} w(t)/\psi(t) = 0$. The function $\psi(s)=cs^\alpha$ with $\alpha>0$, i.e. $\phi(s)= cs^{1+\alpha}$, is an example for a proper growth function.
\end{EX}

\paragraph{Bregman distance.} In order to introduce the notion of a Bregman function \cite{Bregman67}, we first define a set of properties for functions to generate nicely behaving Bregman functions. 
\begin{DEF}[{Legendre function \cite[Def. 5.2]{BBC01}}] 
  The  proper, closed, convex function $\map{\leg}{\X}{\eR}$ is
  \begin{itemize}
    \item[\ii1] \emphdef{essentially smooth}, if $\partial \leg$ is both locally bounded and single-valued on its domain,
    \item[\ii2] \emphdef{essentially strictly convex}, if $(\partial \leg)^{-1}$ is locally bounded on its domain and $\leg$ is strictly convex on every convex subset of $\dom \partial \leg$, and
    \item[\ii3] \emphdef{Legendre}, if $\leg$ is both essentially smooth and essentially strictly convex.
\end{itemize}
\end{DEF}
Note that we have the duality $(\partial \leg)^{-1} = \partial \leg^*$ where $\leg^*$ denotes the conjugate of $\leg$. 
\begin{DEF}[{Bregman distance \cite{Bregman67}}] \label{def:Bregman-distance}
Let $\map{\leg}{\X}{\eR}$ be proper, closed, convex and G\^ateaux differentiable on $\sint\dom\leg \neq \emptyset$. The \emphdef{Bregman distance} associated with $\leg$ is the function 
\[
  \map{\bregmap[\leg]}{\X\times\X}{[0,+\infty]}\,,\qquad 
  (x,\ak{x}) \mapsto \begin{cases}
    \leg(x) - \leg(\ak{x}) - \scal{x-\ak{x}}{\nabla \leg(\ak{x})}\,,&\ \text{if}\ \ak{x}\in \sint\dom\leg \,; \\
    +\infty\,,&\ \text{otherwise}\,.
  \end{cases}
\]
\end{DEF}
In contrast to the Euclidean distance, the Bregman distance is lacking symmetry. 

We focus on Bregman distances that are generated by Legendre functions from the following class:
\[
  \clLeg := \left\{ \map{\leg}{\X}{\eR} \Biggsetsep \matf{h\ \text{is a proper, closed, convex} \\ \text{Legendre function that is}\\ \text{Fr\'echet differentiable on}\ \sint \dom \leg\neq\emptyset } \right\} \,.
\]
To control the variable choice of Bregman distances throughout the algorithm's iterations, we introduce the following ordering relation for $\leg_1,\leg\in\clLeg$:
%we use a fixed Legendre function $\leg\in\clLeg$ for which we define the following relation for $\leg_1,\leg_2\in\clLeg$:
\[
  \leg_1 \succeq \leg \quad \Leftrightarrow \quad \forall x\in \dom \leg\colon \forall \bar x \in \sint\dom\leg\colon\ \breg[\leg_1]{x}{\bar x} \geq \breg[\leg]{x}{\bar x} \,.
\]
As a consequence of $\leg_1\succeq \leg$, we have $\dom\bregmap[\leg_1] \subset \dom\bregmap[\leg]$.
%Moreover, the following subset of $\clLeg$ is introduced:
%\[
%  \clLeg(\leg) := \set{\tilde \leg\in \clLeg \,\vert\, \leg \preceq \tilde \leg } \,.
%\]

In order to conveniently work with Bregman distances, we collect a few properties.
\begin{PROP} \label{prop:basic-Breg-properties}
  Let $\leg\in\clLeg$ and $\bregmap[\leg]$ be the associate Bregman distance.
  \begin{itemize}
    \item[\ii1] $\bregmap[\leg]$ is strictly convex on every convex subset of $\dom\partial \leg$ with respect the first argument.
    \item[\ii2] For $\ak{x}\in \sint\dom \leg$, it holds that $\breg[\leg]x{\ak{x}} = 0$ if and only if $x=\ak{x}$.
    \item[\ii3] For $x\in \X$ and $\bar x, \hat x\in \sint\dom\leg$ the following \emphdef{three point identity} holds:
    \[
      \breg[\leg]{x}{\bar x} = \breg[\leg]{x}{\hat x} + \breg[\leg]{\hat x}{\bar x} + \scal{x - \hat x}{\nabla\leg(\hat x) - \nabla\leg(\bar x)} \,.
    \]
  \end{itemize}
\end{PROP}
\begin{proof}
\ii1 and \ii2 follow directly from the definition of $h$ being essentially strictly convex. \ii3 is stated in \cite{CT93}. It follows from the definition of a Bregman distance.
\end{proof}
Associated with such a distance function is the following proximal mapping.
\begin{DEF}[{Bregman proximal mapping \cite[Def. 3.16]{BBC03}}] \label{def:Bregman-prox}
  Let $\map {f}{\X}{\eR}$ and $\bregmap[\leg]$ be a Bregman distance associated with $\leg\in\clLeg$. The \emphdef{$\bregmap[\leg]$-prox} (or Bregman proximal mapping) associated with $f$ is defined by
  \begin{equation} \label{eq:def-Breg-Prox}
    \Prox[\leg]_{f} (\ak{x}) := \argmin_{x}\, f(x) + \breg[\leg]{x}{\ak{x}} \,.
  \end{equation}
\end{DEF}
In general, the proximal mapping is set-valued, however for a convex function, the following lemma simplifies the situation.
\begin{LEM} \label{lem:breg-prox-single-valued}
  Let $\map{f}{\X}{\eR}$ be a proper, closed, convex function that is bounded from below, and $\leg\in\clLeg$ such that $\sint\dom\leg \cap \dom f\neq\emptyset$. Then the associated Bregman proximal mapping $\Prox[\leg]_{f}$ is single-valued on its domain and maps to $\sint \dom \leg \cap \dom f$.
\end{LEM}
\begin{proof}
Single-valuedness follows from \cite[Corollary~3.25(i)]{BBC03}. The second claim is from \cite[Prop. 3.23(v)(b)]{BBC03}.
\end{proof}
\begin{PROP} \label{prop:three-point-ineq}
  Let $\map{f}{\X}{\eR}$ be a proper, closed, convex function that is bounded from below, and $\leg\in\clLeg$ such that $\sint\dom\leg \cap \dom f\neq\emptyset$. For $\ak{x}\in \sint\dom\leg$, $\hat{x} = \Prox[\leg]_f(\ak{x})$, and any $x\in \dom f$ the following inequality holds:
  \[
      f(x) + \breg[\leg]{x}{\bar x} \geq f(\hat x) + \breg[\leg]{\hat x}{ \bar x} + \breg[\leg]{x}{\hat x} \,.
  \]
\end{PROP}
\begin{proof}
See \cite[Lemma~3.2]{CT93}.
\end{proof}
For examples and more useful properties of Bregman functions, we refer the reader to \cite{BB97,BBC03,BBT16,Nguyen17}. 

\paragraph{Miscellaneous.} We make use of little-$\landau$ notation $f\in\landau(g)$ (or $f=\landau(g)$), which indicates that the asymptotic behavior of a function $f$ is dominated by that of the function $g$. Formally, it is defined by
\[
  f \in \landau(g) \quad \Leftrightarrow\quad   \forall \eps>0\colon \ \abs{f(x)} \leq \eps \abs{g(x)} \ \text{for}\ \vnorm{x}\ \text{sufficiently small}.
\]
Note that a function $\omega$ is in $\o(t)$ if, and only if $\omega$ is a growth function.

% *******************
% >>>>> SECTION <<<<<
% *******************
\section{Line Search Based Bregman Minimization Algorithms} \label{sec:Bregman-alg-ana}

In this paper, we solve optimization problems of the form
\begin{equation} \label{eq:general-problem}
  \min_{x\in\X}\, \fun(x) 
\end{equation}
where $\map{\fun}{\X}{\eR}$ is a proper, closed function on $\X$. We assume that $\Argmin\,f\neq\emptyset$ and $\infbd := \min\fun >-\infty$. The main goal is to develop a provably (subsequentially) convergent algorithm that finds a stationary point $x$ of \eqref{eq:general-problem} in the sense of the limiting slope $\sslp f(x) = 0$.

We analyze abstract algorithms that sequentially minimize convex models of the objective function. 
%The algorithm and analysis is inspired by \cite{BLPP16}. Their analysis is a special case. 

% **********************
% >>>>> SUBSECTION <<<<<
% **********************
\subsection{The Abstract Algorithm} \label{sec:abstract-alg}

%\paragraph{Approximation with local convex models.} 
For each point $\ak{x}$, we consider a proper, closed, convex \emphdef{model function} $\map{\mfun{\ak{x}}}{\X}{\eR}$ with %the property
\begin{equation}\label{eq:model-assumption}
  \abs{\fun(x) - \mfun{\ak{x}}(x)} \leq \omega(\vnorm{x-\ak{x}}) \,,
\end{equation}
where $\omega$ is a growth function as defined in Definition~\ref{def:growth-function}. The \emphdef{model assumption} \eqref{eq:model-assumption} is an abstract description of a (local) first order oracle. For examples, we refer to Section~\ref{sec:examples}. \\

Before delving further into the details, we need a bit of notation. Let
\[
  \mfunp\leg{\ak{x},\ak{z}}(x) := \mfun{\ak{x}}(x) + \breg[\leg]{x}{\ak{z}} 
  \quadand \mfunp\leg{\ak{x}} := \mfunp\leg{\ak{x},\ak{x}} \,,
\]
where $\leg\in\clLeg$. Note that $\mfunp\leg{\ak{x}}(\ak{x})=\fun(\ak{x})$. Moreover, the following quantity defined for generic points $\ak{x}$, $x$ and $\akk{x}$ will be important:
\begin{equation} \label{eq:def-delta}
  \mDeltap\leg{\ak{x}}(x,\akk{x}) := \mfunp\leg{\ak{x}}(x) - \mfunp\leg{\ak{x}}(\akk{x}) \,.
\end{equation}
For $\akk{x}=\ak{x}$, it measures the decrease of the surrogate function $\mfunp\leg{\ak{x}}$ from the current iterate $\ak{x}$ to any other point $x$. Obviously, the definition implies that $\mDeltap\leg{\ak{x}}(x,x) = 0$ for all $x$.

\paragraph{Algorithm.} We consider the following Algorithm~\ref{alg:abstract-alg}. 
\mybox{%
\begin{ALG}[Inexact Bregman Proximal Minimization Line Search] \label{alg:abstract-alg} \
\begin{itemize}
  \item \emph{Basic prerequisites:} Fix $\gamma \in ]0,1[$ and $\leg\in \clLeg$. Let 
  \begin{itemize}
    \item $\seq[\k\in\N]{x\iter\k}$ and $\seq[\k\in\N]{\akk{y}\iter\k}$ be sequences in $\X$;
    \item $\seq[\k\in\N]{\mfun{x\iter\k}}$ be a sequence of model functions with $\inf_{\k\in\N}\inf_{x} \mfun{x\iter\k}(x) > -\infty$;
    \item $\seq[\k\in\N]{\leg\iter\k}$ be a sequence in $\clLeg$ with $\leg\iter\k\succeq \leg$; 
    \item $\seq[\k\in\N]{\eta\iter\k}$ be a sequence of positive real numbers.
  \end{itemize}
  \item \emph{Initialization:} Select $x\iter0 \in \dom\fun\cap \sint\dom\leg$.
  \item \emph{For each $\k\geq 0$:} Generate the sequences such that the following relations hold:
  \begin{eqnarray}
      \label{eq:alg:abstract-relations-A}
      &\mDeltap{\leg\iter\k}{x\iter\k}(\akk{y}\iter\k, x\iter\k) < 0  \ \text{with}\ \akk{y}\iter\k\in \sint\dom\leg\\[1.5ex]
      \label{eq:alg:abstract-relations-B}
      &x\iter\kp = x\iter\k + \eta\iter\k (\akk{y}\iter\k - x\iter\k)  \in \sint\dom \leg \\[1.5ex]
      \label{eq:alg:abstract-relations-C}
      &\fun(x\iter\kp) \leq \fun(x\iter\k) + \gamma \eta\iter\k \mDeltap{\leg\iter\k}{x\iter\k}(\akk{y}\iter\k, x\iter\k)
  \end{eqnarray}
  If \eqref{eq:alg:abstract-relations-A} cannot be satisfied, then the algorithm terminates.
\end{itemize}
\end{ALG}
}
The algorithm starts with a feasible point\footnote{It is often easy to find a feasible point. Of course, there are cases, where finding an initialization is a problem itself. We assume that the user provides a feasible initial point.} $x\iter0$. At each iteration, it computes a point $\akk{y}\iter\k$ that satisfies \eqref{eq:alg:abstract-relations-A}, which is an inexact solution of the \emphdef{Bregman proximal mapping}
\begin{equation}\label{eq:Breg-prox-map}
  \akk{x}\iter\k=\Prox[\leg\iter\k]_{\mfun{x\iter\k}}(x\iter\k) := \argmin_{x\in \X} \mfun{x\iter\k}(x) + \breg[\leg\iter\k]{x}{x\iter\k}
\end{equation}
that, at least, improves the (model) value compared to $x\iter\k$. Thanks to the class of Legendre functions $\clLeg$, this proximal mapping is well-defined and single-valued on its domain. The exact version of the algorithm solves the proximal mapping exactly for the global optimal solution. The optimal solution of the proximal mapping will always be denoted by $\akk{x}\iter\k$ instead of $\akk{y}\iter\k$, which refers to an approximate solution. The direction $\akk{y}\iter\k-x\iter\k$ can be considered as a descent direction for the function $\fun$. Given this direction, the goal of \eqref{eq:alg:abstract-relations-B} and \eqref{eq:alg:abstract-relations-C} is the estimation of a step size $\eta\iter\k$ (by line search, cf. Algorithm~\ref{alg:line-search}) that reduces the value of the objective function. In case that the proximal mapping has a solution but the first relation \eqref{eq:alg:abstract-relations-A} can only be satisfied with equality, we will see that $x\iter\k=\akk{y}\iter\k$ must be a stationary point of the objective, hence, the algorithm terminates. 
\begin{REM}  \label{rem:alg-backtracking-subproblems}
Instead of performing backtracking on the objective values as in Algorithm~\ref{alg:abstract-alg}, backtracking on the scaling of the Bregman distance in \eqref{eq:alg:abstract-relations-A} is also possible. For a special model function, this leads to ProxDescent \cite{LW16,DL16} (with Euclidean proximity function). If a scaled version of \eqref{eq:alg:abstract-relations-A} yields a descent on $\fun$, we can set $\eta\iter\k=1$, and accept this point. However, this can be expensive when the proximal subproblem in \eqref{eq:alg:abstract-relations-A} is hard to solve, since each trial step requires to solve the subproblem. In order to break the backtracking, the new objective value must be computed anyway. Therefore, a computational advantage of the line search \eqref{eq:alg:abstract-relations-B} and \eqref{eq:alg:abstract-relations-C} is to be expected (cf. Section~\ref{sec:numeric-robust-regression}). 
\end{REM}
\mybox{%
\begin{ALG}[Line Search for Algorithm~\ref{alg:abstract-alg}] \label{alg:line-search} \ 
\begin{itemize}
  \item \emph{Basic prerequisites:} Fix $\delta,\gamma \in ]0,1[$,  $\tilde\eta>0$, and $\k\in \N$. 
  \item \emph{Input:} Current iterates $x\iter\k\in\sint\dom\leg$ and $\akk{y}\iter\k$ satisfy \eqref{eq:alg:abstract-relations-A}.
  \item \emph{Solve:} Find the smallest $j\in \N$ such that $\tilde\eta\iter{j}:=\tilde\eta\delta^j$ satisfies \eqref{eq:alg:abstract-relations-B} and \eqref{eq:alg:abstract-relations-C}.
  \item \emph{Return:} Set the feasible step size $\eta\iter\k$ for iteration $\k$ to $\tilde\eta\iter{j}$.
\end{itemize}
\end{ALG}
}

Algorithm~\ref{alg:abstract-alg}--\ref{alg:line-search} is well defined as the following lemmas show.
\begin{LEM}[Well-definedness] \label{lem:alg-well-defined} 
  Let $\omega$ in \eqref{eq:model-assumption} be a growth function. Algorithm~\ref{alg:abstract-alg} is well-defined, i.e., for all $\k\in \N$, the following holds:
  \begin{itemize}
    \item[\ii1] there exists $\akk{y}\iter\k$ that satisfies \eqref{eq:alg:abstract-relations-A} or $x\iter\k=\akk{x}\iter\k$ and the algorithm terminates; 
    \item[\ii2] $x\iter\k \in \dom \fun \cap \sint\dom\leg$; and 
    \item[\ii3] there exists $\eta\iter\k$ that satisfies \eqref{eq:alg:abstract-relations-B} and \eqref{eq:alg:abstract-relations-C}.
  \end{itemize}
\end{LEM}
\begin{proof}
  \ii1 For $x\iter\k \in \sint \dom \leg$, Lemma~\ref{lem:breg-prox-single-valued} shows that $\Prox[\leg\iter\k]_{\mfun{x\iter\k}}$ maps to $\sint\dom\leg\iter\k\cap \dom \mfun{x\iter\k}\subset\sint\dom\leg\cap \dom \fun$ and is single-valued. Thus, for example, $\akk{y}\iter\k = \akk{x}\iter\k$ satisfies \eqref{eq:alg:abstract-relations-A}. Otherwise, $x\iter\k = \akk{x}\iter\k$, which shows \ii1.
  \ii2 Since $x\iter0\in \dom\fun\cap\sint\dom\leg$ and $\fun(x\iter\kp) \leq \fun(x\iter\k)$ by \eqref{eq:alg:abstract-relations-C} it holds that $x\iter\k \in \dom \fun$ for all $\k$. Since $x\iter\k\in\sint\dom\leg$ and $\akk{y}\iter\k\in \dom\leg$, for small $\eta\iter\k$ also $x\iter\kp\in \sint\dom\leg$, hence $x\iter\kp \in \dom \fun\cap\sint\dom\leg$. Inductively, we conclude the statement.
  \ii3 This will be shown in Lemma~\ref{lem:finite-escape}.
\end{proof}
\begin{LEM}[Finite termination of Algorithm~\ref{alg:line-search}] \label{lem:finite-escape}
  Consider Algorithm~\ref{alg:abstract-alg} and fix $\k\in\N$. Let $\omega$ in \eqref{eq:model-assumption} be a growth function. Let $\delta, \gamma\in ]0,1[$, $\tilde\eta>0$, $\ak{\leg}:=\leg\iter\k$, and $\ak{x}:=x\iter\k$, $\akk{y}:=\akk{y}\iter\k$ be such that $\mDeltap{\ak{\leg}}{\ak{x}}(\akk{y},\ak{x}) < 0$. Then, there exists $j\in\N$ such that $\tilde\eta\iter{j}:=\tilde\eta\delta^j$ satisfies  
  \[
      f(\ak{x} + \tilde\eta\iter{j} (\akk{y}-\ak{x})) \leq f(\ak{x}) + \gamma\tilde\eta\iter{j} \mDeltap{\ak{\leg}}{\ak{x}}(\akk{y},\ak{x}) \,.
  \]
\end{LEM}
\begin{proof}
This result is proved by contradiction. Define $\dir:=\akk{y}-\ak{x} $. By our assumption in \eqref{eq:model-assumption}, we observe that
\begin{equation} \label{eq:prop:finite-exit-of-Armijo-A}
  f(\ak{x}+\tilde\eta\iter{j} \dir) - f(\ak{x}) \leq \mfun{\ak{x}}(\ak{x}+\tilde\eta\iter{j} \dir) - f(\ak{x}) + \landau (\tilde\eta\iter{j}) \,.
\end{equation}
Using Jensen's inequality for the convex function $\mfun{\ak{x}}$ provides:
\begin{equation} \label{eq:prop:finite-exit-of-Armijo-B}
  \mfun{\ak{x}}(\ak{x} + \tilde\eta\iter{j} \dir) -\mfun{\ak{x}}(\ak{x})\leq \tilde\eta\iter{j} \mfun{\ak{x}}(\ak{x} + \dir) + (1-\tilde\eta\iter{j})\mfun{\ak{x}}(\ak{x}) -\mfun{\ak{x}}(\ak{x}) = \tilde\eta\iter{j}\cdot\Big( \mfun{\ak{x}}(\ak{x} + \dir) - \mfun{\ak{x}}(\ak{x}) \Big) \,.
\end{equation}
Now, suppose $\gamma \mDeltap{\ak{\leg}}{\ak{x}}(\akk{y},\ak{x}) < \frac{1}{\tilde\eta\iter{j}} ( f(\ak{x} + \tilde\eta\iter{j} \dir) - f(\ak{x}) )$  holds for any $j\in\N$. Then, using \eqref{eq:prop:finite-exit-of-Armijo-A} and \eqref{eq:prop:finite-exit-of-Armijo-B}, we conclude the following: 
\[
  \begin{split}
    \gamma \mDeltap{\ak{\leg}}{\ak{x}}(\akk{y},\ak{x}) 
    %<&\ \frac{1}{\tilde\eta\iter{j}} ( f(\ak{x} + \tilde\eta\iter{j} \dir) - f(\ak{x}) ) \\
    <&\ \mfun{\ak{x}}(\ak{x} + \dir) - \mfun{\ak{x}}(\ak{x}) + \landau(\tilde\eta\iter{j})/ \tilde\eta\iter{j} \\
    \leq&\ \mfun{\ak{x}}(\ak{x} + \dir) - \mfun{\ak{x}}(\ak{x}) + \breg[{\ak{\leg}}]{\akk{y}}{\ak{x}} + \landau(\tilde\eta\iter{j})/ \tilde\eta\iter{j}\\
    =&\ (\mfunp{\ak{\leg}}{ \ak{x}}(\akk{y})-\mfunp{\ak{\leg}}{\ak{x}}(\ak{x})) + \landau(\tilde\eta\iter{j})/ \tilde\eta\iter{j} = \mDeltap{\ak{\leg}}{\ak{x}}(\akk{y},\ak{x})+\landau(\tilde\eta\iter{j})/ \tilde\eta\iter{j} \,,
  \end{split}
\]
which for $j\to\infty$ yields the desired contradiction, since $\gamma \in ]0,1[$ and $\mDeltap{\ak{\leg}}{\ak{x}}(\akk{y},\ak{x})<0$.
\end{proof}
%Let us now analyze Algorithm~\ref{alg:abstract-alg}.
%Before, we analyze Algorithm~\ref{alg:abstract-alg} in Section~\ref{sec:convergence-ana-finite}, we provide examples.

% **********************
% >>>>> SUBSECTION <<<<<
% **********************
\subsection{Finite Time Convergence Analysis} \label{sec:convergence-ana-finite}

First, we study the case when the algorithm terminates after a finite number of iterations, i.e., there exists $\k_0\in\N$ such that \eqref{eq:alg:abstract-relations-A} cannot be satisfied. Then, the point $\akk{y}\iter{\k_0}$ is a global minimizer of $\mfunp{\leg\iter{\k_0}}{x\iter{\k_0}}$ and $\mDeltap{\leg\iter{\k_0}}{x\iter{\k_0}}(\akk{y}\iter{\k_0}, x\iter{\k_0}) =0$. Moreover, the point $x\iter{\k_0}$ turns out to be a stationary point of $\fun$. 
\begin{LEM} \label{lem:subdiff-equal} 
  For $\ak{x}\in\dom \fun$ and a model $\mfun{\ak{x}}$ that satisfies \eqref{eq:model-assumption}, where $\omega$ is a growth function, the following holds:
  \[
        \slp{\mfun{\ak{x}}}(\ak{x}) = \slp{\fun}(\ak{x}) \,.
  \]
\end{LEM}
\begin{proof} Since $\omega(0)=0$, we have from \eqref{eq:model-assumption} that $\mfun{\ak{x}}(\ak{x}) = \fun(\ak{x})$. This, together with sub-additivity of $\pos{\cdot}$, entails
\begin{align*}
\frac{\pos{\mfun{\ak{x}}(\ak{x}) - \mfun{\ak{x}}(x)}}{\vnorm{x-\ak{x}}} \leq \frac{\pos{\fun(\ak{x}) - \fun(x)} + \pos{\fun(x) - \mfun{\ak{x}}(x)}}{\vnorm{x-\ak{x}}} 
&\leq \frac{\pos{\fun(\ak{x}) - \fun(x)}}{\vnorm{x-\ak{x}}} + \frac{\abs{\fun(x) - \mfun{\ak{x}}(x)}}{\vnorm{x-\ak{x}}} \\
&\leq \frac{\pos{\fun(\ak{x}) - \fun(x)}}{\vnorm{x-\ak{x}}} + \frac{\omega(\vnorm{x-\ak{x}})}{\vnorm{x-\ak{x}}}
\end{align*}
Passing to the $\limsup$ on both sides and using that $\omega \in \landau(t)$, we get
\[
\slp{\mfun{\ak{x}}}(\ak{x}) \leq \slp{\fun}(\ak{x}) .
\]
Arguing similarly but now starting with $\slp{\fun}(\ak{x})$, we get the reverse inequality, which in turn shows the claimed equality.
%\[
%\slp{\mfun{\ak{x}}}(\ak{x}) =  \limsup_{\substack{x\to \ak{x}\\x\neq\ak{x}}} \frac{\pos{\mfun{\ak{x}}(\ak{x}) - \mfun{\ak{x}}(x)}}{\vnorm{x-\ak{x}}} = \limsup_{\substack{x\to\ak{x}\\x\neq\ak{x}}} \frac{\pos{\fun(\ak{x}) - \fun(x)} + \pos{\fun(x) - \mfun{\ak{x}}(x)}}{\vnorm{x-\ak{x}}} = \slp{\fun}(\ak{x}) \,,
%\]
%since $\pos{\fun(x) - \mfun{\ak{x}}(x)}\leq \landau(\vnorm{x-\ak{x}})$.
\end{proof}
\begin{PROP}[Stationarity for finite time termination] \label{prop:equivalence-stationarity}
  Consider the setting of Algorithm~\ref{alg:abstract-alg}. Let $\omega$ in \eqref{eq:model-assumption} be a growth function. Let $\k_0\in\N$ be fixed, and set $\akk{x}=\akk{y}\iter{\k_0}$, $\ak{x} = x\iter{\k_0}$, $\ak{\leg}=\leg\iter{\k_0}$, and $\ak{x},\akk{x}\in \dom\fun\cap \sint\dom\leg$. If $\mDeltap{\ak{\leg}}{\ak{x}}(\akk{x},\ak{x})\geq 0$, then $\akk{x}=\ak{x}$, $\mDeltap{\ak{\leg}}{\ak{x}}(\akk{x},\ak{x}) =  0$, and $\slp \fun(\ak{x})=0$, i.e. $\ak{x}$ is a stationary point of $\fun$.
\end{PROP}
\begin{proof}
  Since $\akk{x}$ is the unique solution of the proximal mapping, obviously $\mDeltap{\ak{\leg}}{\ak{x}}(\akk{x},\ak{x}) =  0$ and $\akk{x}=\ak{x}$. Moreover, $\akk{x}$ is the minimizer of $\mfunp{\ak{\leg}}{\ak{x}}$, i.e. we have
\[
  0 = \slp{\mfunp{\ak{\leg}}{\ak{x}}}(\akk{x}) = \slp{\mfunp{\ak{\leg}}{\ak{x}}}(\ak{x}) = \limsup_{\substack{x\to\ak{x}\\x\neq\ak{x}}} \frac{\pos{\mfun{\ak{x}}(\ak{x}) - \mfun{\ak{x}}(x) - \breg[{\ak{\leg}}]{x}{\ak{x}} }}{\vnorm{x-\ak{x}}} = \slp{\mfun{\ak{x}}}(\ak{x}) = \slp{\fun}(\ak{x}) \,,
\]
where we used that $\ak{\leg}$ is Fr\'echet differentiable at $\ak{x}$ and Lemma~\ref{lem:subdiff-equal}. 
\end{proof}

% **********************
% >>>>> SUBSECTION <<<<<
% **********************
\subsection{Asymptotic Convergence Analysis} \label{sec:convergence-ana-asymptotic}

 We have established stationarity of the algorithm's output, when it terminates after a finite number of iterations. Therefore, without loss of generality, we now focus on the case where \eqref{eq:alg:abstract-relations-A} can be satisfied for all $\k\in\N$. We need to make the following assumptions.
\begin{ASS} \label{ass:vanishing-prox-val-error}
  The sequence $\seq[\k\in\N]{\akk{y}\iter\k}$ satisfies $\mfunp{\leg\iter\k}{x\iter\k}(\akk{y}\iter\k) \leq \inf \mfunp{\leg\iter\k}{x\iter\k} + \eps\iter\k$ for some $\eps\iter\k\to0$.
\end{ASS}
\begin{REM}
  Assumption~\ref{ass:vanishing-prox-val-error} states that asymptotically (for $\k\to\infty$) the Bregman proximal mapping \eqref{eq:Breg-prox-map} must be solved accurately. In order to obtain stationarity of a limit point, Assumption~\ref{ass:vanishing-prox-val-error} is necessary, as shown by Bonettini et al. \cite[after Theorem~4.1]{BLPP16} for a special setting of model functions.
\end{REM}
\begin{ASS} \label{ass:Bregman-continuity}
  Let $\leg\in\clLeg$. For every bounded sequences $\seq[\k\in\N]{x\iter\k}$ and $\seq[\k\in\N]{\ak{x}\iter\k}$ in $\sint\dom\leg$, and $\seq[\k\in\N]{\leg\iter\k}$ such that $\leg\iter\k\succeq\leg$, it is assumed that:
  \[
      x\iter\k - \ak{x}\iter\k \to 0 \quad\Leftrightarrow\quad \breg[\leg\iter\k]{x\iter\k}{\ak{x}\iter\k}\to 0 \,.
  \]
\end{ASS}
\begin{REM}
  \begin{itemize}
  \item[\ii1] Assumption~\ref{ass:Bregman-continuity} states that (asymptotically) a vanishing Bregman distance reflects a vanishing Euclidean distance. This is a natural assumption and satisfied, e.g., by most entropies such as Boltzmann--Shannon, Fermi--Dirac, and Burg entropy. 
  \item[\ii2] The equivalence in Assumption~\ref{ass:Bregman-continuity} is satisfied, for example, when there exists $c\in \R$ such that $c\,\leg \succeq \leg\iter\k$ holds for all $\k\in\N$ and the following holds:
  \[
      x\iter\k - \ak{x}\iter\k \to 0 \quad\Leftrightarrow\quad \breg[\leg]{x\iter\k}{\ak{x}\iter\k}\to 0 \,.
  \]
  \end{itemize}
\end{REM}

%Convergence of the objective values is immediate by construction as we now summarize.
\begin{PROP}[Convergence of objective values]  \label{prop:value-convergence}
  Consider the setting of Algorithm~\ref{alg:abstract-alg}. Let $\omega$ in \eqref{eq:model-assumption} be a growth function. The sequence of objective values $\seq[\k\in\N]{\fun(x\iter\k)}$ is non-increasing and converging to some $\fun^*\geq \infbd>-\infty$.
\end{PROP}
\begin{proof}
  This statement is a consequence of \eqref{eq:alg:abstract-relations-C} and \eqref{eq:alg:abstract-relations-A}, and the lower-boundedness of $\fun$.
\end{proof}

 Asymptotically, under some condition on the step size, the improvement of the model objective value between $\akk{y}\iter\k$ and $x\iter\k$ must tend to zero. Since we do not assume that the step sizes $\eta\iter\k$ are bounded away from zero, this is a non-trivial result.
\begin{PROP}[Vanishing model improvement]  \label{prop:model-convergence}
  Consider the setting of Algorithm~\ref{alg:abstract-alg}. Let $\omega$ in \eqref{eq:model-assumption} be a growth function. Suppose, either\footnote{Note that $\inf_\k\eta\iter\k>0$ is equivalent to $\liminf_\k\eta\iter\k>0$ as we assume $\eta\iter\k>0$ for all $\k\in\N$.} $\inf_\k\eta\iter\k>0$ or $\eta\iter\k$ is selected by the Line Search Algorithm~\ref{alg:line-search}. Then, 
  \[
    \ds\sum_{\k=0}^\infty  \eta\iter\k (-\mDeltap{\leg\iter\k}{x\iter\k}(\akk{y}\iter\k, x\iter\k)) < +\infty
    \quad \text{and}\quad  \mDeltap{\leg\iter\k}{x\iter\k}(\akk{y}\iter\k, x\iter\k) \to 0 \ \text{ as $\k\to\infty$}.
  \]
\end{PROP}
\begin{proof}
 The first part follows from rearranging \eqref{eq:alg:abstract-relations-C}, and summing both sides for $\k=0,\ldots,n$:
  \[
    \gamma \sum_{\k=0}^n  \eta\iter\k (-\mDeltap{\leg\iter\k}{x\iter\k}(\akk{y}\iter\k, x\iter\k)) \leq \sum_{\k=0}^n (\fun(x\iter\k) - \fun(x\iter\kp) ) = \fun(x\iter0) - \fun(x\iter{n+1}) \leq \fun(x\iter0) - \fun^* \,.
  \]
  In the remainder of the proof, we show that $\mDeltap{\leg\iter\k}{x\iter\k}(\akk{y}\iter\k, x\iter\k) \to 0$, which is not obvious unless $\inf_\k\eta\iter\k>0$. The model improvement is bounded. Boundedness from above is satisfied by construction of the sequence $\seq[\k\in\N]{\akk{y}^\k}$. Boundedness from below follows from the following observation and the uniform boundedness of the model functions from below:
\[
    \mDeltap{\leg\iter\k}{x\iter\k}(\akk{y}\iter\k, x\iter\k) 
    = \mfunp{\leg\iter\k}{{x}\iter\k}(\akk{y}\iter\k) - \mfunp{\leg\iter\k}{{x}\iter\k}({x}\iter\k) 
    \geq  \mfunp{\leg\iter\k}{{x}\iter\k}(\akk{x}\iter\k) - \fun({x}\iter\k)
    \geq  \mfun{{x}\iter\k}(\akk{x}\iter\k) - \fun({x}\iter0)\,.
\]
Therefore, there exists $K\subset\N$ such that the subsequence $\mDeltap{\leg\iter\k}{x\iter\k}(\akk{y}\iter\k, x\iter\k)$ converges to some $\Delta^*$ as $\k\fto {K}\infty$. Suppose $\Delta^*<0$. Then, the first part of the statement implies that the step size sequence must tend to zero, i.e., $\eta\iter\k\to 0$ for $\k\fto{K}\infty$. For $\k\in K$ sufficiently large, the line search procedure in Algorithm~\ref{alg:line-search} reduces the step length from $\eta\iter\k/\delta$ to $\eta\iter\k$. (Note that $\eta\iter\k$ can be assumed to be the \enquote{first} step length that achieves a reduction in \eqref{eq:alg:abstract-relations-C}). Before multiplying with $\delta$, no descent of \eqref{eq:alg:abstract-relations-C} was observed, i.e.,
\[
  ({\eta\iter\k}/\delta) \gamma \mDeltap{\leg\iter\k}{x\iter\k}(\akk{y}\iter\k, x\iter\k) < \fun(x\iter\k + (\eta\iter\k/\delta) \dir\iter\k) - \fun(x\iter\k)\,,
\]
where $\dir\iter\k = \akk{y}\iter\k - x\iter\k$. Using \eqref{eq:prop:finite-exit-of-Armijo-A} and \eqref{eq:prop:finite-exit-of-Armijo-B}, we can make the same observation as in the proof of Lemma~\ref{lem:finite-escape}:
\[
  \begin{split}
    \gamma \mDeltap{\leg\iter\k}{x\iter\k}(\akk{y}\iter\k,x\iter\k) 
    <&\ \mfun{x\iter\k}(x\iter\k + \dir) - \mfun{x\iter\k}(x\iter\k) + \landau(\eta\iter\k/\delta)/ (\eta\iter\k/\delta) \\
    \leq&\ \mfun{x\iter\k}(x\iter\k + \dir) - \mfun{x\iter\k}(x\iter\k) + \breg[{\leg\iter\k}]{\akk{y}\iter\k}{x\iter\k} + \landau(\eta\iter\k)/ \eta\iter\k\\
    =&\ (\mfunp{\leg\iter\k}{ x\iter\k}(\akk{y}\iter\k)-\mfunp{\leg\iter\k}{x\iter\k}(x\iter\k)) + \landau(\eta\iter\k)/ \eta\iter\k \\
    =&\  \mDeltap{\leg\iter\k}{x\iter\k}(\akk{y}\iter\k,x\iter\k)+\landau(\eta\iter\k)/ \eta\iter\k \,,
  \end{split}
\]
which for $\eta\iter\k \to 0$ yields a contradiction, since $\gamma \in ]0,1[$ and $\mDeltap{\leg\iter\k}{x\iter\k}(\akk{y}\iter\k,x\iter\k)<0$. Therefore, any cluster point $\Delta^*$ of $\seq[\k\in K]{\mDeltap{\leg\iter\k}{x\iter\k}(\akk{y}\iter\k, x\iter\k)}$ must be $0$, which concludes the proof.
\end{proof}

\subsubsection{Asymptotic Stationarity with a Growth Function} \label{sec:conv-ana-growth-fun}

In order to establish stationarity of limit points generated by Algorithm~\ref{alg:abstract-alg} additional assumptions are required. We consider three different settings for the model assumption \eqref{eq:model-assumption}: $\omega$ in the model assumption \eqref{eq:model-assumption} is a growth function (this section), $\omega$ is a proper growth function (Section~\ref{sec:conv-ana-proper-growth-fun}), and $\omega$ is global growth function of the form $\omega=\bregmap[\leg]$ (Section~\ref{sec:global-model}).
\begin{ASS}\label{ass:vanishing-slope}
Let $x^*$ be a limit point of $\seq[\k\in\N]{x\iter\k}$ and $x\iter\k \fto{\fun} x^*$ as $\k\fto K\infty$ with $K\subset\N$. Then 
\[
  \slp{\mfun{x\iter\k}}(x\iter\k) = \slp{\fun}(x\iter\k) \to 0 \quad \text{as}\ \k\fto K\infty \,.
\]
\end{ASS}
\begin{REM} \label{rem:vanishing-slope} 
  Assumption~\ref{ass:vanishing-slope} is common for abstract algorithms. Attouch et al. \cite{ABS13}, for example, use a relative error condition of the form $\slp{\fun}(x\iter\kp) \leq b\vnorm{x\iter\kp - x\iter\k}$, $b\in\R$. A weaker sufficient condition for Assumption~\ref{ass:vanishing-slope} is $\slp{\fun}(x\iter\kp) \leq \psi(\vnorm{x\iter\kp - x\iter\k})$ for some continuous function $\map{\psi}{\R_+}{\R_+}$ with $\psi(0)=0$; See Corollary~\ref{cor:stationarity-growth}. See also Remark~\ref{rem:alt-to-ass:vanishing-slope}. For explicit examples, we refer to Section~\ref{subsec:examples-models}.
\end{REM}
Using this assumption, we can state one of our main theorems, which shows convergence to a stationary point under various condition. The conditions are easily verified in many applications (see Section~\ref{sec:examples}). 
\begin{THM}[Asymptotic stationarity with a growth function] \label{thm:stationarity-growth}
  Consider the setting of Algorithm~\ref{alg:abstract-alg}. Let $\omega$ in \eqref{eq:model-assumption} be a growth function. Moreover, let either $\inf_\k\eta\iter\k>0$ or $\eta\iter\k$ be selected by the Line Search Algorithm~\ref{alg:line-search}. Let $\seq[\k\in\N]{x\iter\k}$ and $\seq[\k\in\N]{\akk{y}\iter\k}$ be bounded sequences such that Assumptions~\ref{ass:vanishing-prox-val-error} and~\ref{ass:Bregman-continuity} hold and let $\mfun{x\iter\k}$ obey \eqref{eq:model-assumption} with growth function $\omega$. Then, $x\iter\k - \akk{y}\iter\k \to 0$ and for $\akk{x}\iter\k=\Prox[\leg\iter\k]_{\mfun{x\iter\k}}(x\iter\k)$, it holds that $x\iter\k - \akk{x}\iter\k \to 0$ and $\akk{x}\iter\k - \akk{y}\iter\k \to 0$. Moreover, $\fun(x\iter\k) - \fun(\akk{y}\iter\k)\to 0$ and $\fun(\akk{x}\iter\k)-\fun(x\iter\k) \to 0$ as $\k\to\infty$. Suppose Assumption~\ref{ass:vanishing-slope} is satisfied. If $x^*$ is a limit point of the sequence $\seq[\k\in\N]{x\iter\k}$, and one of the following conditions is satisfied:
  \begin{itemize}
  \item[\ii1] $\fun$ is continuous on the closure of $\dom\leg$,
  \item[\ii2] $x^*\in\sint\dom\leg$,
  \item[\ii3] $x^* \in \dom\leg$ and $\breg[\leg\iter\k]{x^*}{\akk{y}\iter\k} \to 0$ as $\k\fto{K}\infty$,
  \item[\ii4] $x^*\in \scl\dom\leg$ and 
  \begin{itemize}
    \item for all $x\in \sint\dom\leg\cap \dom\fun$ holds that $\breg[\leg\iter\k]{x}{\akk{x}\iter\k}-\breg[\leg\iter\k]{x}{x\iter\k} \to 0$ as $\k\fto{K}\infty$,
    \item and for all $x\in \dom\fun$ the model functions obey $\mfun{x\iter\k}(x)\to \mfun{x^*}(x)$ as $\k\fto{K}\infty$,
  \end{itemize}
  \end{itemize}
  then $x^*$ is a stationary point of $\fun$.
\end{THM}  
\begin{proof}
  First, we show that for $\k\to \infty$ the pairwise distances between the sequences $\seq[\k\in\N]{x\iter\k}$, $\seq[\k\in\N]{\akk{y}\iter\k}$, and $\seq[\k\in\N]{\akk{x}\iter\k}$ vanishes. Proposition~\ref{prop:three-point-ineq}, reformulated in our notation, can be stated as
  \begin{equation} \label{eq-B:thm:stationarity-growth}
    \mDeltap{\leg\iter\k}{x\iter\k}(x,\akk{x}\iter\k) = \mfunp{\leg\iter\k}{x\iter\k}(x) - \mfunp{\leg\iter\k}{x\iter\k}(\akk{x}\iter\k) \geq \breg[\leg\iter\k]{x}{\akk{x}\iter\k}\,, \quad \forall x\in \dom \fun \,. %\geq \breg[\leg]{x}{\akk{x}\iter\k}\,.
  \end{equation}
  As a direct consequence, using $x=\akk{y}\iter\k$ together with Assumptions~\ref{ass:vanishing-prox-val-error} and~\ref{ass:Bregman-continuity}, we obtain
  \[%\begin{equation} \label{eq-C:thm:stationarity-growth}
    \breg[\leg\iter\k]{\akk{y}\iter\k}{\akk{x}\iter\k}\to0 \quad\text{thus}\quad \akk{x}\iter\k - \akk{y}\iter\k \to 0 \,.
  \]%\end{equation}
  Moreover, from Proposition~\ref{prop:model-convergence}, we have $\mDeltap{\leg\iter\k}{x\iter\k}(\akk{y}\iter\k, x\iter\k) \to 0$, and from
  \begin{equation} \label{eq-D:thm:stationarity-growth}
    \mDeltap{\leg\iter\k}{x\iter\k}(\akk{y}\iter\k,x\iter\k) = \mDeltap{\leg\iter\k}{x\iter\k}(\akk{y}\iter\k,\akk{x}\iter\k) - \mDeltap{\leg\iter\k}{x\iter\k}(x\iter\k,\akk{x}\iter\k) \leq \mDeltap{\leg\iter\k}{x\iter\k}(\akk{y}\iter\k,\akk{x}\iter\k)-\breg[\leg\iter\k]{x\iter\k}{\akk{x}\iter\k} \,,
  \end{equation}
  and Assumptions~\ref{ass:vanishing-prox-val-error} and~\ref{ass:Bregman-continuity}, we conclude that $x\iter\k - \akk{x}\iter\k \to 0$, hence also $x\iter\k - \akk{y}\iter\k \to 0$. \\

  The next step is to show that $\fun(x\iter\k) - \fun(\akk{y}\iter\k)\to 0$ as $k\to \infty$. This follows from the following estimation:
  \begin{equation} \label{eq-E:thm:stationarity-growth}
    \begin{split}
    \abs{\fun(x\iter\k) - \fun(\akk{y}\iter\k)} \leq&\  \abs{\mfun{x\iter\k}(x\iter\k) - \mfun{x\iter\k}(\akk{y}\iter\k)} + \omega(\vnorm{\akk{y}\iter\k-x\iter\k}) \\
    \leq&\  \abs{\mfunp{\leg\iter\k}{x\iter\k}(x\iter\k) - \mfunp{\leg\iter\k}{x\iter\k}(\akk{y}\iter\k)} + \breg[\leg\iter\k]{\akk{y}\iter\k}{x\iter\k} + \omega(\vnorm{\akk{y}\iter\k-x\iter\k}) \\
    =&\  \abs{\mDeltap{\leg\iter\k}{x\iter\k}(x\iter\k,\akk{y}\iter\k)} + \breg[\leg\iter\k]{\akk{y}\iter\k}{x\iter\k} + \omega(\vnorm{\akk{y}\iter\k-x\iter\k}) \,,
    \end{split}
  \end{equation}
  where the right hand side vanishes for $\k\to\infty$. Analogously, we can show that $\fun(\akk{x}\iter\k)-\fun(x\iter\k) \to 0$ as $\k\to\infty$. \\

%%  The next step is to show that $\fun(x\iter\k) - \fun(\akk{y}\iter\k)\to 0$ as $k\to \infty$. This follows from the following estimations, which use the model assumption \eqref{eq:model-assumption}
%%  \[
%%    \begin{split}
%%    \fun(\akk{y}\iter\k) - \fun(x\iter\k) \leq&\  \mfun{x\iter\k}(\akk{y}\iter\k) - \mfun{x\iter\k}(x\iter\k) + \omega(\vnorm{\akk{y}\iter\k-x\iter\k}) \\
%%    =&\  \mfunp{\leg\iter\k}{x\iter\k}(\akk{y}\iter\k) - \mfunp{\leg\iter\k}{x\iter\k}(x\iter\k) - \breg[\leg\iter\k]{\akk{y}\iter\k}{x\iter\k} + \omega(\vnorm{\akk{y}\iter\k-x\iter\k}) \\
%%    \leq&\  \mDeltap{\leg\iter\k}{x\iter\k}(\akk{y}\iter\k,x\iter\k) + \omega(\vnorm{\akk{y}\iter\k-x\iter\k}) \\
%%    \end{split}
%%  \]
%%  where the right hand side vanishes for $\k\to\infty$. %Analogously, we can show that $\fun(\akk{x}\iter\k)-\fun(x\iter\k) \to 0$ as $\k\to\infty$. \\

  Let $x^*$ be the limit point of the subsequence $\seq[\k\in K]{x\iter\k}$ for some $K\subset \N$. The remainder of the proof shows that $\fun(\akk{y}\iter\k)\to f(x^*)$ as $\k\to\infty$. Then $\fun(x\iter\k) - \fun(\akk{y}\iter\k)\to 0$ implies that $x\iter\k\fto{\fun} x^*$ as $\k\fto K\infty$, and by Assumption~\ref{ass:vanishing-slope}, the slope vanishes, hence the limiting slope $\sslp{\fun}(x^*)$ at $x^*$ also vanishes, which concludes the proof.\\

  \ii1 implies $\fun(\akk{y}\iter\k)\to f(x^*)$ as $\k\to\infty$ by definition. For \ii2 and \ii3, we make the following observation:
  \begin{equation} \label{eq-A:thm:stationarity-growth}
    \fun(\akk{y}\iter\k) - \omega(\vnorm{\akk{y}\iter\k-x\iter\k}) \leq \mfunp{\leg\iter\k}{x\iter\k}(\akk{y}\iter\k) = \mfunp{\leg\iter\k}{x\iter\k}(\akk{x}\iter\k) + (\mfunp{\leg\iter\k}{x\iter\k}(\akk{y}\iter\k) - \mfunp{\leg\iter\k}{x\iter\k}(\akk{x}\iter\k)) \leq \mfunp{\leg\iter\k}{x\iter\k}(x^*) + \mDeltap{\leg\iter\k}{x\iter\k}(\akk{y}\iter\k,\akk{x}\iter\k) \,,
  \end{equation}
  where $\akk{x}\iter\k = \Prox[\leg\iter\k]_{\mfun{x\iter\k}}(x\iter\k)$. Taking \enquote{$\limsup_{\k\fto K\infty}$} on both sides, $\breg[\leg\iter\k]{x^*}{x\iter\k} \to 0$ (Assumption~\ref{ass:Bregman-continuity} for \ii2 or the assumption in \ii3), and $\mDeltap{\leg\iter\k}{x\iter\k}(\akk{y}\iter\k,\akk{x}\iter\k)\to 0$ (Assumption~\ref{ass:vanishing-prox-val-error}) shows that $\limsup_{\k\fto K\infty} \fun(\akk{y}\iter\k) \leq \fun(x^*)$. Since $\fun$ is closed, $\fun(\akk{y}\iter\k) \to \fun(x^*)$ holds. 

  We consider \ii4. For all $x\in \sint\dom\leg\cap \dom\fun$, we have \eqref{eq-B:thm:stationarity-growth} or, reformulated, $\mfunp{\leg\iter\k}{x\iter\k}(x) - \breg[\leg\iter\k]{x}{\akk{x}\iter\k}  \geq  \mfunp{\leg\iter\k}{x\iter\k}(\akk{x}\iter\k)$, which implies the following:
  \[
        \mfun{x\iter\k}(x) + \breg[\leg\iter\k]{x}{x\iter\k} - \breg[\leg\iter\k]{x}{\akk{x}\iter\k} - \breg[\leg\iter\k]{\akk{x}\iter\k}{x\iter\k}  \geq  \fun(\akk{x}\iter\k) - \omega(\vnorm{\akk{x}\iter\k- x\iter\k})   \,.
  \]
  Note that for any $x$ the limits for $k\fto{K}\infty$ on the left hand side exist. In particular, we have 
  \[
    \breg[\leg\iter\k]{x}{x\iter\k} - \breg[\leg\iter\k]{x}{\akk{x}\iter\k} - \breg[\leg\iter\k]{\akk{x}\iter\k}{x\iter\k} \to 0\quad \text{as}\ \k\fto{K}\infty \,,
  \]
  by the assumption in \ii4, and Assumption~\ref{ass:Bregman-continuity} together with $\akk{x}\iter\k-x\iter\k\to 0$. The limit of $\mfun{x\iter\k}(x)$ exists by assumption and coincides with $\mfun{x^*}(x)$. Choosing a sequence $\seq[\k\in\N]{z\iter\k}$ in $\sint\dom\leg\cap \dom\fun$ with $z\iter\k \to x^*$ as $k\fto{K}\infty$, in the limit, we obtain
  \[
    \fun(x^*) \geq \lim_{\k\fto{K}\infty} \fun(\akk{x}\iter\k) =: \fun^*\,,
  \]
  since $\mfun{x^*}(z\iter\k) \to \mfun{x^*}(x^*) = \fun(x^*)$ for $z\iter\k \to x^*$ as $\k\fto{K}\infty$. Invoking that $\fun$ is closed, we conclude the $f$-attentive convergence $\fun^* = \liminf_{\k\to\infty} \fun(x\iter\k) \geq  \fun(x^*)\geq \fun^*$.
\end{proof}
\begin{REM}
Existence of a limit point $x^*$ is guaranteed by assuming that $\seq[\k\in\N]{x\iter\k}$ is bounded. Alternatively, we could require that $\fun$ is coercive (i.e. $\fun(x\iter\k)\to \infty$ for $\vnorm{x\iter\k}\to\infty$), which implies boundedness of the lower level sets of $\fun$, hence by Proposition~\ref{prop:value-convergence} the boundedness of $\seq[\k\in\N]{x\iter\k}$.
\end{REM}
\begin{REM} \label{rem:alt-to-ass:vanishing-slope}
From Theorem~\ref{thm:stationarity-growth}, clearly, also $\akk{y}\iter\k \fto\fun x^*$ and $\akk{x}\iter\k\fto\fun x^*$ as $\k\fto{K}\infty$ holds. Therefore, Assumption~\ref{ass:vanishing-slope} could also be stated as the requirement
\[
  \slp{\fun}(\akk{x}\iter\k) \to 0 \quad \text{or}\quad \slp{\fun}(\akk{y}\iter\k)\to0\quad \text{as}\ \k\fto K\infty \,,
\]
in order to conclude that limit points of $\seq[\k\in\N]{x\iter\k}$ are stationary points.
\end{REM}
As a simple corollary of this theorem, we replace Assumption~\ref{ass:vanishing-slope} with the relative error condition mentioned in Remark~\ref{rem:vanishing-slope}.
\begin{COR}[Asymptotic stationarity with a growth function] \label{cor:stationarity-growth}
  Consider the setting of Algorithm~\ref{alg:abstract-alg}. Let $\omega$ in \eqref{eq:model-assumption} be a growth function. Moreover, let either $\inf_\k\eta\iter\k>0$ or $\eta\iter\k$ be selected by the Line Search Algorithm~\ref{alg:line-search}. Let $\seq[\k\in\N]{x\iter\k}$ and $\seq[\k\in\N]{\akk{y}\iter\k}$ be bounded sequences such that Assumptions~\ref{ass:vanishing-prox-val-error} and~\ref{ass:Bregman-continuity} hold and let $\mfun{x\iter\k}$ obey \eqref{eq:model-assumption} with growth function $\omega$. Suppose there exists a continuous function $\map{\psi}{\R_+}{\R_+}$ with $\psi(0)=0$ such that $\slp{\fun}(x\iter\kp) \leq \psi(\vnorm{x\iter\kp - x\iter\k})$ is satisfied. If $x^*$ is a limit point of the sequence $\seq[\k\in\N]{x\iter\k}$ and one of the conditions \ii1--\ii4 in Theorem~\ref{thm:stationarity-growth} is satisfied, then $x^*$ is a stationary point of $\fun$.
\end{COR}  
\begin{proof}
  Theorem~\ref{thm:stationarity-growth} shows that $\akk{y}\iter\k - x\iter\k \to 0$, thus, $\inf_\k \eta\iter\k>0$ implies $x\iter\kp-x\iter\k \to 0$ by \eqref{eq:alg:abstract-relations-B}. Therefore, the relation $\slp{\fun}(x\iter\kp) \leq \psi(\vnorm{x\iter\kp - x\iter\k})$ shows that Assumption~\ref{ass:vanishing-slope} is automatically satisfied and we can apply Theorem~\ref{thm:stationarity-growth} to deduce the statement. 
\end{proof}

\paragraph{Some more results on the limit point set.} In Theorem~\ref{thm:stationarity-growth} we have shown that limit points of the sequence $\seq[\k\in\N]{x\iter\k}$ generated by Algorithm~\ref{alg:abstract-alg} are stationary, and in fact the sequence $\fun$-converges to its limit points. The following proposition shows some more properties of the set of limit points of $\seq[\k\in\N]{x\iter\k}$. This is a well-known result \cite[Lem.~5]{BST14} that follows from $x\iter\kp - x\iter\k \to 0$ as $\k\to\infty$.
\begin{PROP} \label{prop:limit-point-set-connected}
  Consider the setting of Algorithm~\ref{alg:abstract-alg}. Let $\omega$ in \eqref{eq:model-assumption} be a growth function and $\inf_\k\eta\iter\k>0$. Let $\seq[\k\in\N]{x\iter\k}$ and $\seq[\k\in\N]{\akk{y}\iter\k}$ be bounded sequences such that Assumptions~\ref{ass:vanishing-prox-val-error}, \ref{ass:Bregman-continuity} and~\ref{ass:vanishing-slope} hold. Suppose one of the conditions \ii1--\ii4 in Theorem~\ref{thm:stationarity-growth} is satisfied for each limit point of $\seq[\k\in\N]{x\iter\k}$. Then, the set 
  $%\[
    \mathfrak S := \set{ x^*\in\X\setsep \exists K\subset\N\colon  x\iter\k \to  x^*\ \text{as}\ \k\fto K\infty }
  $ %\]
  of limit points of $\seq[\k\in\N]{x\iter\k}$
  is connected, each point $x^*\in \mathfrak S$ is stationary for $\fun$, and $\fun$ is constant on $\mathfrak S$.
\end{PROP}
\begin{proof}
  Theorem~\ref{thm:stationarity-growth} shows that $\akk{y}\iter\k - x\iter\k \to 0$. Thus, boundedness of $\eta\iter\k$ away from $0$ implies $x\iter\kp-x\iter\k \to 0$ by \eqref{eq:alg:abstract-relations-B}. Now, the statement follows from \cite[Lem.~5]{BST14} and Theorem~\ref{thm:stationarity-growth}.
\end{proof}

\subsubsection{Asymptotic Stationarity with  a Proper Growth Function} \label{sec:conv-ana-proper-growth-fun}

 Our proof of stationarity of limit points generated by Algorithm~\ref{alg:abstract-alg} under the assumption of a proper growth function $\omega$ in \eqref{eq:model-assumption} relies on an adaptation of a recently proved result by Drusvyatskiy et al. \cite[Corollary~5.3]{DIL16}, which is stated in Lemma~\ref{cor:general-stationarity} before the main theorem of this subsection. The credits for this lemma should go to \cite{DIL16}.
\begin{LEM}[Perturbation result under approximate optimality] \label{cor:general-stationarity}
 Let $\map{\fun}{\X}{\eR}$ be a proper closed function. Consider bounded sequences $\seq[\k\in\N]{x\iter\k}$ and $\seq[\k\in\N]{\akk{y}\iter\k}$ with $x\iter\k-\akk{y}\iter\k\to 0$ for $\k\to\infty$, and model functions $\mfunp{\leg\iter\k}{x\iter\k}$ according to \eqref{eq:model-assumption} with proper growth functions. Suppose Assumption~\ref{ass:vanishing-prox-val-error} and~\ref{ass:Bregman-continuity} hold. If $(x^*,\fun(x^*))$ is a limit point of $\seq[\k\in\N]{x\iter\k,\fun(x\iter\k)}$, then $x^*$ is stationary for $\fun$.
\end{LEM}
\begin{proof}
Recall $\eps\iter\k$ from Assumption~\ref{ass:vanishing-prox-val-error}. Theorem 5.1 from \cite{DIL16} guarantees, for each $\k$ and any $\rho\iter\k>0$, the existence of points $\akp{y}\iter\k$ and $z\iter\k$ such that the following hold:
  \begin{itemize}
    \item[\ii1] (point proximity)
    \[
      \vnorm{\akk{y}\iter\k-z\iter\k} \leq \frac{\eps\iter\k}{\rho\iter\k}\quad\text{and}\quad \vnorm{z\iter\k - \akp{y}\iter\k} \leq 2\cdot \frac{\omega(\vnorm{z\iter\k-x\iter\k})}{\omega_+^\prime(\vnorm{z\iter\k-x\iter\k})}\,,
    \]
    under the convention $\frac 00=0$\,,
    \item[\ii2] (value proximity) $\fun(\akp{y}\iter\k) \leq \fun(\akk{y}\iter\k) + 2\omega(\vnorm{z\iter\k-x\iter\k}) + \omega(\vnorm{\akk{y}\iter\k-x\iter\k})$, and
    \item[\ii3] (near-stationarity) $\slp\fun(\akp{y}) \leq \rho\iter\k + \omega_+^\prime(\vnorm{z\iter\k - x\iter\k}) + \omega_+^\prime(\vnorm{\akp{y}\iter\k-x\iter\k})$,
  \end{itemize}
  Setting $\rho\iter\k=\sqrt{\eps\iter\k}$, using $\eps\iter\k\to0$ and the point proximity, shows that $\vnorm{\akk{y}\iter\k-z\iter\k}\to0$. Moreover $\vnorm{z\iter\k - x\iter\k} \leq \vnorm{z\iter\k-\akk{y}\iter\k} + \vnorm{\akk{y}\iter\k - x\iter\k}\to 0$, which implies that $\vnorm{z\iter\k-\akp{y}\iter\k}\to 0$. Now, we fix a convergent subsequence $(x\iter\k,\fun(x\iter\k)) \to (x^*,\fun(x^*))$ as $\k\fto K\infty$ for some $K\subset\N$. Using \eqref{eq-B:thm:stationarity-growth}, we observe $\akk{x}\iter\k-\akk{y}\iter\k\to 0$, hence $x\iter\k-\akk{x}\iter\k\to 0$. From Proposition~\ref{prop:model-convergence} and Assumption~\ref{ass:vanishing-prox-val-error}, we conclude that $\mDeltap{\leg\iter\k}{x\iter\k}(\akk{y}\iter\k,x\iter\k)\to 0$, and, therefore $\fun(x\iter\k)-\fun(\akk{y}\iter\k)\to0$ using \eqref{eq-E:thm:stationarity-growth}. Consider the value proximity. Combined with the lower semi-continuity of $\fun$, it yields
  \[
    \fun(x^*) \leq \liminf_{\k\fto K\infty} \fun(\akp{y}\iter\k) \leq \limsup_{\k\fto K\infty} \fun(\akp{y}\iter\k) \leq \limsup_{\k\fto K\infty} \fun(\akk{y}\iter\k) \leq f(x^*) \,,
  \]
  hence $(\akp{y}\iter\k,\fun(\akp{y}\iter\k)) \to (x^*,\fun(x^*))$ as $\k\fto K\infty$. Near-stationarity implies that $\slp\fun(\akp{y}\iter\k) \to 0$, which proves that $\sslp\fun(x^*)=0$, hence $x^*$ is a stationary point.
\end{proof}
\begin{REM}
  The setting in \cite[Corollary~5.3]{DIL16} is recovered when $\seq[\k\in\N]{x\iter\k}$ is given by $x\iter\kp = \akk{y}\iter\k$.
\end{REM}
\begin{THM}[Asymptotic stationarity with a proper growth function] \label{thm:stationarity-proper-growth}
  Consider the setting of Algorithm~\ref{alg:abstract-alg}. Let $\omega$ in \eqref{eq:model-assumption} be a proper growth function. Moreover, let either $\inf_\k\eta\iter\k>0$ or $\eta\iter\k$ be selected by the Line Search Algorithm~\ref{alg:line-search}. Let $\seq[\k\in\N]{x\iter\k}$ and $\seq[\k\in\N]{\akk{y}\iter\k}$ be bounded sequences such that Assumptions~\ref{ass:vanishing-prox-val-error} and~\ref{ass:Bregman-continuity} hold. If $x^*$ is a limit point of the sequence $x\iter\k$ and one of the conditions \ii1--\ii4 in Theorem~\ref{thm:stationarity-growth} is satisfied, then $x^*$ is a stationary point of $\fun$.
\end{THM}  
\begin{proof}
Propositions~\ref{prop:value-convergence} and~\ref{prop:model-convergence}, and the proof of $\fun$-attentive convergence from Theorem~\ref{thm:stationarity-growth} only rely on a growth function. Instead of assuming that the slope vanishes, here we apply Lemma~\ref{cor:general-stationarity} to conclude stationarity of the limit points.
\end{proof}

Of course, Proposition~\ref{prop:limit-point-set-connected} can also be stated in the context here.

\subsubsection{Asymptotic Analysis with a Global Growth Function} \label{sec:global-model}

Suppose, for $\ak{x}\in\sint\dom\leg$ for some $\leg\in \clLeg$, the model error can be estimated as follows:
\begin{equation}\label{eq:model-assumption-global}
  \abs{\fun(x) - \mfun{\ak{x}}(x)} \leq L \breg[{\leg}]{x}{\ak{x}} \quad \forall x \,.
\end{equation}
Since $\leg$ is Fr\'echet differentiable on $\sint\dom\leg$, the right hand side is bounded by\xspace a growth function. Without loss of generality, we restrict ourselves to a fixed function $\leg\in\clLeg$ (this section analyses a single iteration). In order to reveal similarities to well-known step size rules, we scale $\leg$ in the definition of $\mfunp{\leg}{\ak{x}}$ to $\bregmap[{\leg}/\alpha]=\frac{1}{\alpha}\bregmap[{\leg}]$ with $\alpha>0$ instead of $\bregmap[{\leg}]$. Here, decreasing objective values can be assured without the line search procedure (see Proposition~\ref{prop:global:descent-property}), i.e., $\eta\iter\k = 1$ is always feasible. 

In order to obtain the result of stationarity of limit points (Theorem~\ref{thm:stationarity-growth} or \ref{thm:stationarity-proper-growth}), we can either verify by hand that Assumption~\ref{ass:vanishing-slope} holds or we need to assume that $\breg[\leg]{x}{\ak{x}}$ is bounded by\xspace a proper growth function. 
\begin{PROP} \label{prop:global:descent-property}
  Consider the setting of Algorithm~\ref{alg:abstract-alg} and let \eqref{eq:model-assumption-global} be satisfied. 
  \begin{itemize}
    \item[\ii1] For points $\akk{y}$ that satisfy $\mDeltap{\leg}{\ak{x}}(\akk{y},\ak{x}) < 0$, 
  \[
     \frac {1-\alpha L}{\alpha} \breg[{\leg}]{\akk{y}}{\ak{x}} \leq \fun(\ak{x}) - \fun(\akk{y}) 
     %({1-\alpha L})/{\alpha}\cdot \breg[{\leg}]{\akk{y}}{\ak{x}} \leq \fun(\ak{x}) - \fun(\akk{y}) 
  \]
  holds, where the left-hand-side is strictly larger than $0$ for $\alpha\in ]0,1/L[$.
  \item[\ii2] For points $\akk{x}=\Prox[\leg]_{\mfun{\ak{x}}}(\ak{x})$, the following descent property holds:
  \[
     \frac{1+\rho-\alpha L}\alpha \breg[{\leg}]{\akk{x}}{\ak{x}} \leq \fun(\ak{x}) - \fun(\akk{x}) \,,
     %({1+\rho-\alpha L})/\alpha\cdot \breg[{\leg}]{\akk{x}}{\ak{x}} \leq \fun(\ak{x}) - \fun(\akk{x}) \,,
  \]
  where the left-hand-side is strictly larger than $0$ for $\alpha\in ]0,(1+\rho)/L[$, and $\rho$ is the Bregman symmetry factor defined by 
  $\rho := \inf\set{\frac{\breg[\leg] x{\bar x}}{\breg[\leg]{\bar x}{x}} \setsep x,\bar x \in \sint\dom \leg\,, \ x\neq \bar x }$;  (see \cite{BBT16}).
 \end{itemize}
\end{PROP}
\begin{proof}
The following relations hold:
\begin{equation} \label{eq:prop:global:descent-property-a}
  \mDeltap{\leg}{\ak{x}}(\akk{y},\ak{x}) \leq 0
  \quad\Leftrightarrow\quad \mfunp{\leg}{\ak{x}}(\akk{y}) \leq \mfunp{\leg}{\ak{x}}(\ak{x}) 
  \quad \Leftrightarrow \quad \mfun{\ak{x}}(\akk{y}) + \frac 1\alpha\breg[{\leg}]{\akk{y}}{\ak{x}} \leq \mfun{\ak{x}}(\ak{x}) = \fun(\ak{x}) \,.
\end{equation}
Bounding the left hand side of the last expression using \eqref{eq:model-assumption-global}, we obtain
\begin{equation} \label{eq:prop:global:descent-property-b}
   \fun(\akk{y}) - L\breg[{\leg}]{\akk{y}}{\ak{x}} + \frac1\alpha\breg[{\leg}]{\akk{y}}{\ak{x}} \leq \fun(\ak{x}) \,,
\end{equation}
which proves part \ii1. Part \ii2 follows analogously. However, thanks to the three point inequality from Proposition~\ref{prop:three-point-ineq} and optimality of $\akk{x}$ the rightmost inequality of \eqref{eq:prop:global:descent-property-a} improves to
\[
  \mfun{\ak{x}}(\akk{x}) + \frac 1\alpha\breg[{\leg}]{\akk{x}}{\ak{x}} + \frac 1\alpha\breg[{\leg}]{\ak{x}}{\akk{x}} \leq \mfun{\ak{x}}(\ak{x}) = \fun(\ak{x}) \,,
\]
and the statement follows.
\end{proof}

% *******************
% >>>>> SECTION <<<<<
% *******************
\subsection{A Remark on Convex Optimization} 

In this section, let $\fun$ be convex, and consider the following global model assumption
\begin{equation}\label{eq:model-assumption-conv}
  0 \leq \fun(x) - \mfun{\ak{x}}(x) \leq L\breg[\leg]{x}{\ak{x}} \,.
\end{equation}
We establish a convergence rate of $\O(1/\k)$ for Algorithm~\ref{alg:abstract-alg} with $\eta\iter\k\equiv1$. For Forward--Backward Splitting, this has been shown by Bauschke et al. \cite{BBT16}. We only require $\mfun{\ak{x}}$ to be a model w.r.t. \eqref{eq:model-assumption-conv}. 
%If we restrict ourselves to $\leg$ being the squared Euclidean distance, we can accelerated the algorithm to show a rate $\O(1/\k^2)$. 
\begin{PROP} \label{prop:conv-rate-convex}
  Consider Algorithm~\ref{alg:abstract-alg} with $\eta\iter\k\equiv1$ and model functions that obey \eqref{eq:model-assumption-conv}. For $x\iter\kp=\Prox[\leg/\alpha]_{\mfun{x\iter\k}}(x\iter\k)$ and $\alpha = \frac 1L$, the following rate of convergence on the objective values holds:
 \[
   \fun(x\iter\kp) - \fun(x) \leq \frac{L \breg[\leg]{x^*}{x\iter0}}{2 \k}   \quad (=\O(1/\k))\,.
 \]
\end{PROP}
\begin{proof}
The three point inequality in Proposition~\ref{prop:three-point-ineq} combined with the model assumption \eqref{eq:model-assumption-conv} yields the following inequality:
\[
  \fun(\akk{x}) + \frac{1-\alpha L}\alpha \breg[\leg]{\akk{x}}{\ak{x}} + \frac 1\alpha \breg[\leg]{x}{\akk{x}}  \leq \fun(x) + \frac{1}{\alpha} \breg[\leg]{x}{\ak{x}} 
\]
for all $x$. Restricting to $0<\alpha\leq \frac 1L$, we obtain
\begin{equation}\label{eq:descent-relation}
  \fun(\akk{x}) - \fun(x) \leq \frac{1}{\alpha} \left( \breg[\leg]{x}{\ak{x}} - \breg[\leg]{x}{\akk{x}}  \right) \,.
\end{equation}
Let $x^*$ be a minimizer of $\fun$. We make the following choices:
\[
  x=x^*\,,\quad \akk{x}=x\iter\kp\,, \quadand \ak{x} = x\iter\k \,.
\]
Summing both sides up to iteration $\k$ and the descent property yield the convergence rate:
\begin{equation}\label{eq:abstract-conv}
  \fun(x\iter\kp) - \fun(x) \leq \frac{\breg[\leg]{x^*}{x\iter0}}{2\alpha \k } \overset{\alpha=\frac 1L}{=} \frac{L \breg[\leg]{x^*}{x\iter0}}{2 \k }\,.
\end{equation}
\end{proof}
\begin{EX}
This section goes beyond the Forward--Backward Splitting setting, e.g., we may specialize Example~\ref{ex:prox-descent} to convex problems of the form $\min_{x} g(F(x))$, for instance, using $F(x)  = (f_1(x),f_2(x))$ and $g(z_1,z_2)=\max\{z_1,z_2\}$.
\end{EX}

% *******************
% >>>>> SECTION <<<<<
% *******************
\section{Examples} \label{sec:examples}

We discuss several classes of problems that can be solved using our framework. To apply Algorithm~\ref{alg:abstract-alg}, in Section~\ref{subsec:examples-models}, we define a suitable model and mention the associated algorithmic step that arises from exactly minimizing the sum of the model and an Euclidean proximity measure. However our algorithm allows for inexact solutions and very flexible (also iteration dependent) Bregman proximity functions. Examples are provided in Section~\ref{subsec:examples-bregman}. For brevity, we define the symbols \emphdef{$\Gamma_0$ for the set of proper, closed, convex functions} and \emphdef{$C^1$ for the set of continuously differentiable functions}.

\subsection{Examples of Model Functions} \label{subsec:examples-models}

\begin{EX}[Forward--Backward Splitting] \label{ex:fbs}
Problems of the form 
\[
  \fun=f_0 + f_1\qquad \text{with}\ f_0\in\Gamma_0\ \text{and}\ f_1\in C^1
\]
can be modeled by 
\[
  \mfun{\ak{x}}(x)=f_0(x) + f_1(\ak{x}) + \scal{x-\ak{x}}{\nabla f_1(\ak{x})}\,.
\]
This model is associated with Forward--Backward Splitting (FBS). We assume that one of the following error models is satisfied:
\[
  \abs{\fun(x) - \mfun{\ak{x}}(x)} = \abs{f_1(x) - f_1(\ak{x}) - \scal{x-\ak{x}}{\nabla f_1(\ak{x})}}
  \leq \begin{cases}
    \frac L2\vnorm{x-\ak{x}}^2\,,& \text{if} \ \nabla f_1\ \text{is L-Lipschitz}\,; \\[1mm]
    \int_0^1\frac{\phi(t\vnorm{x-\ak{x}})}{t}\,\mathit{d}t  \,,& 
    \begin{minipage}{3.5cm}if  $\nabla f_1$ is $\psi$-uniform continuous\,;\end{minipage} \\[1mm]
    \omega(\vnorm{x-\ak{x}})\,, & \text{otherwise}\,,
  \end{cases}
\]
which is the linearization error of the smooth part $f_1$. The first case obeys a global (proper) growth function, derived from the common Descent Lemma. The second case is the generalization to a $\psi$-uniformly continuous gradient with $\phi(s)=s\psi(s)$ as in Lemma~\ref{lem:gen-descent-lemma}. The bound $\int_0^1\frac{\phi(t\vnorm{x-\ak{x}})}{t}\,\mathit{d}t$ is a growth function but not necessarily a proper growth function. The third case is the most general and assumes that the error obeys a growth function. In any case, the model satisfies the model consistency required in Theorem~\ref{thm:stationarity-growth}\ii4. For any $x\in \dom\fun$ and $\ak{x} \to x^*$, 
\[
  \abs{f_0(x) + f_1(\ak{x}) + \scal{x-\ak{x}}{\nabla f_1(\ak{x})} - (f_0(x) + f_1(x^*) + \scal{x-x^*}{\nabla f_1(x^*)} )} \to 0
\]
holds, thanks to the continuous differentiability of $f_1$ and continuity of the inner product.\\

In order to verify Assumption~\ref{ass:vanishing-slope}, we make use of Remark~\ref{rem:alt-to-ass:vanishing-slope} and show that $\slp{\fun}(\akk{x}\iter\k) \to 0$ as $\k\fto{K}\infty$ where $K\subset \N$ is such that $x\iter\k\fto{K}x^*$. Note that $\akk{x}\iter\k$ satisfies the following relation:
\begin{gather*}
  0\in  \partial f_0(\akk{x}\iter\k) + \nabla f_1(x\iter\k) + \nabla \leg\iter\k(\akk{x}\iter\k) - \nabla \leg\iter\k(x\iter\k) \\
\Rightarrow \nabla f_1(\akk{x}\iter\k) - \nabla f_1(x\iter\k) + \nabla \leg\iter\k(x\iter\k) - \nabla \leg\iter\k(\akk{x}\iter\k) \in  \partial f_0(\akk{x}\iter\k) +  \nabla f_1(\akk{x}\iter\k) = \partial \fun(\akk{x}\iter\k)
\end{gather*}
Moreover, we know that $\akk{x}\iter\k-x\iter\k\to 0$ as $\k\fto{K}\infty$. Since $\nabla f_1$ is continuous, if $\vnorm{\nabla \leg\iter\k(x\iter\k) - \nabla \leg\iter\k(\akk{x}\iter\k)}\to 0$ for $\k\fto{K}\infty$, then Assumption~\ref{ass:vanishing-slope}/Remark~\ref{rem:alt-to-ass:vanishing-slope} is satisfied. The condition $\vnorm{\nabla \leg\iter\k(x\iter\k) - \nabla \leg\iter\k(\akk{x}\iter\k)}\to 0$ is naturally fulfilled by many Legendre functions, e.g., if $\nabla \leg\iter\k$ is $\psi$-uniformly continuous\xspace (uniformly in $\k$) with $\alpha>0$ or uniformly continuous (independent of $\k$) on bounded sets or continuous at $x^*$ (uniformly w.r.t. $\k$), and will be discussed in more detail in Section~\ref{subsec:examples-bregman}.
\end{EX}
\begin{EX}[Variable metric FBS] \label{ex:vm-fbs}
We consider an extension of Examples~\ref{ex:fbs}. An alternative feasible model for a twice continuously differentiable function $f_1$ is the following:
\[
  \mfun{\ak{x}}(x)=f_0(x) + f_1(\ak{x})+ \scal{x-\ak{x}}{\nabla f_1(\ak{x})}+ \frac 12\scal{x-\ak{x}}{B(x-\ak{x})}\,,
\]
where $B:=[\nabla^2 f_1(\ak{x})]_+$ is a positive definite approximation to $\nabla^2 f_1(\ak{x})$, which leads to a Hessian driven variable metric FBS. It is easy to see that the model error satisfies the growth function $\omega(s)$. Again, Theorem~\ref{thm:stationarity-growth}\ii4 obviously holds and the same conclusions about Assumption~\ref{ass:vanishing-slope} can be made as in Example~\ref{ex:fbs}.
\end{EX}
\begin{EX}[ProxDescent] \label{ex:prox-descent}
Problems of the form 
\[
  f_0 + g\circ F\qquad \text{with}\ f_0\in\Gamma_0\,,\  F\in C^1\,,\ \text{and}\  g\in \Gamma_0\  \text{finite-valued}\,, 
\]
which often arise from non-linear inverse problems, can be approached by the model function
\[
  \mfun{\ak{x}}(x)=f_0(x) + g(F(\ak{x}) +  DF(\ak{x}) (x-\ak{x}))\,,
\]
where $DF(\ak{x})$ is the Jacobian matrix of $F$ at $\ak{x}$. The associated algorithm is connected to ProxDescent \cite{LW16,DL16}.  If $g$ is a quadratic function, the algorithm reduces to the Levenberg--Marquardt algorithm~\cite{Marquardt63}. The error model can be computed as follows:
\begin{equation} \label{eq:ex:ProxDescent-model-error}
  \begin{split}
  \abs{\fun(x) - \mfun{\ak{x}}(x)} 
    =&\  \abs{g(F(x)) - g(F(\ak{x}) + DF(\ak{x}) (x-\ak{x})) } \\
    \leq&\  \ell \vnorm{F(x) - F(\ak{x}) - DF(\ak{x}) (x-\ak{x})} \\
    \leq&\  \begin{cases}
         \frac {\ell L}2\vnorm{x-\ak{x}}^2\,,& \text{if} \ DF\ \text{is L-Lipschitz}\ \text{and}\ g\ \text{is $\ell$-Lipschitz} \,; \\[1mm]
    \ell\int_0^1\frac{\phi(t\vnorm{x-\ak{x}})}{t}\,\mathit{d}t  \,,& 
    \text{if}\ DF\ \text{is $\psi$-uniform continuous}\ \text{and}\ g\ \text{is $\ell$-Lipschitz}\,; \\[1mm]
         \omega(\vnorm{x-\ak{x}})\,, & \text{otherwise}\,,
       \end{cases}
  \end{split}
\end{equation}
where $\ell$ is the (possibly local) Lipschitz constant of $g$ around $F(\ak{x})$. Since $g$ is convex and finite-valued, it is always locally Lipschitz continuous. Since $F$ is continuously differentiable, for $x$ sufficiently close to $\ak{x}$, both $F(x)$ and $F(\ak{x}) + DF(\ak{x}) (x-\ak{x})$ lie in a neighborhood of $F(\ak{x})$ where the local Lipschitz constant $\ell$ of $g$ is valid, which shows the first inequality in \eqref{eq:ex:ProxDescent-model-error}. The last line in \eqref{eq:ex:ProxDescent-model-error} assumes that, either the linearization error of $F$\xspace obeys a global proper growth function, the specific growth function $\int_0^1\frac{\phi(t\vnorm{x-\ak{x}})}{t}\,\mathit{d}t$,\xspace or it obeys a growth function $\omega(s)$. With a similar reasoning, we can show that Theorem~\ref{thm:stationarity-growth}\ii4 is satisfied. \\

We consider Assumption~\ref{ass:vanishing-slope} (see also Remark~\ref{rem:alt-to-ass:vanishing-slope}). Let $x\iter\k \to x^*$ as $\k\fto{K}\infty$ for $K\subset \N$ and $\akk{x}\iter\k-x\iter\k\to 0$. Since $g$ is finite-valued, using \cite[Corollary~16.38]{BC11} (sum-rule for the subdifferential), and \cite[Theorem~10.6]{Rock98}, we observe that
\begin{equation} \label{eq-A:ex:prox-descent}
  0\in \partial f_0(\akk{x}\iter\k) + DF(x\iter\k)^* \partial g(F(x\iter\k) + DF(x\iter\k)(\akk{x}\iter\k-x\iter\k)) + \nabla \leg\iter\k(\akk{x}\iter\k) - \nabla \leg\iter\k(x\iter\k) \,,
\end{equation}
where $DF(x\iter\k)^*$ denotes the adjoint of $DF(x\iter\k)$. We can assume that, for $\k$ large enough, $F(x\iter\k) + DF(x\iter\k)(\akk{x}\iter\k-x\iter\k)$ and $F(\akk{x}\iter\k)$ lie a neighborhood of $F(x^*)$ on which $g$ has the Lipschitz constant $\ell>0$. By \cite[Theorem~9.13]{Rock98}, $\partial g$ is locally bounded around $F(x^*)$, i.e. there exists a compact set $G$ such that $\partial g(z) \subset G$ for all $z$ in a neighborhood of $F(x^*)$. We conclude that
\[
  \sup_{\substack{ v\in \partial g(F(x\iter\k) + DF(x\iter\k)(\akk{x}\iter\k-x\iter\k))\\
                   w\in \partial g(F(\akk{x}\iter\k)) } } \vnorm{ DF(x\iter\k)^* v - DF(\akk{x}\iter\k)^* w} 
  \leq \sup_{v,w\in G} \vnorm{ DF(x\iter\k)^* v - DF(\akk{x}\iter\k)^* w} \to 0
\]
for $\k\fto{K}\infty$  since $DF(x\iter\k)\to DF(x^*)$ and $DF(\akk{x}\iter\k)\to DF(x^*)$. Again assuming that $\nabla \leg\iter\k(\akk{x}\iter\k) - \nabla \leg\iter\k(x\iter\k) \to 0$ we conclude that the outer set-limit of the right hand side of \eqref{eq-A:ex:prox-descent} is included in $\partial f(\akk{x}\iter\k)$ and, therefore, the slope $\slp{\fun}(\akk{x}\iter\k)$ vanishes for $k\fto{K}\infty$.
\end{EX}
\begin{EX} \label{ex:mm}
  Problems of the form 
  \[
    f_0 + g\circ F\qquad \text{with}\ f_0\in\Gamma_0\,,\ g\in C^1
    \,, \nabla g_i \geq 0
    \,,\ \text{and}\ F=(F_1,\ldots,F_M)\  \text{is Lipschitz with}\ F_i\in\Gamma_0
  \]
  can be modeled by 
  \[
    \mfun{\ak{x}}(x) = f_0(x)+g(F(\ak{x})) + \scal{F(x) - F(\ak{x})}{\nabla g(F(\ak{x}))} \,.
  \]
  Such problems appear for example in non-convex regularized imaging problems in the context of iteratively reweighted algorithms \cite{ODBP15}. For the error of this model function, we observe the following:
  \[
    \begin{split}
    \abs{\fun(x) - \mfun{\ak{x}}(x)} 
    =&\  \abs{g(F(x)) - ( g(F(\ak{x})) + \scal{F(x) - F(\ak{x})}{\nabla g(F(\ak{x}))})} \\
    \leq &\  \begin{cases} 
            \frac \ell2 \abs{ F(x) - F(\ak{x}) }^2\,,& \text{if}\ \nabla g\ \text{is $\ell$-Lipschitz} \,; \\
            \int_0^1\frac{\phi(t\vnorm{F(x)-F(\ak{x})})}{t}\,\mathit{d}t\,,& \text{if}\ \nabla g\ \text{is $\psi$-uniform continuous} \,; \\
            \omega(\vnorm{ F(x) -  F(\ak{x}) })\,, & \text{otherwise}\,;  
          \end{cases} \\
    \leq &\  \begin{cases} 
            \frac {\ell L^2}2 \abs{ x - \ak{x} }^2\,,& \text{if}\ \nabla g\ \text{is $\ell$-Lipschitz}\ \text{and}\ F\ \text{is $L$-Lipschitz} \,; \\[1mm]
            \int_0^1\frac{\phi(tL\vnorm{x-\ak{x}}))}{t}\,\mathit{d}t\,, & 
              \begin{minipage}{7.5cm}if $\nabla g$ is $\psi$-uniform continuous and\\ $F$ is $L$-Lipschitz  \,;\end{minipage} \\[1mm]
            \omega(c\vnorm{ x -  \ak{x} })\,, & \text{otherwise}\,,
          \end{cases} \\
    \end{split}
  \]
  which shows the same growth functions are obeyed as in Example~\ref{ex:fbs} and~\ref{ex:prox-descent}. The explanation for the validity of the reformulations are analogue to those of Example~\ref{ex:prox-descent}. It is easy to see that Theorem~\ref{thm:stationarity-growth}\ii4 holds.\\
  
We consider Assumption~\ref{ass:vanishing-slope}/Remark~\ref{rem:alt-to-ass:vanishing-slope}. Let $x\iter\k \to x^*$ as $\k\fto{K}\infty$ for $K\subset\N$ and $\akk{x}\iter\k-x\iter\k\to 0$. Since $g$ is continuously differentiable, the sum-formula for the subdifferential holds. Moreover, we can apply \cite[Corollary~10.09]{Rock98} (addition of functions) to see that $\akk{x}\iter\k$ satisfies the following relation:
\[
  0 \in \partial f_0(\akk{x}\iter\k) + \sum_{i=1}^M \partial F_{i}(\akk{x}\iter\k) (\nabla g(F(x\iter\k)))_i + \nabla \leg\iter\k(\akk{x}\iter\k) - \nabla \leg\iter\k(x\iter\k) \,,
\]
Note that by \cite[Theorem~10.49]{Rock98} the subdifferential of $g \circ F$ at $\akk{x}\iter\k$ is $\sum_{i=1}^M \partial F_{i}(\akk{x}\iter\k) (\nabla g(F(\akk{x}\iter\k)))_i$. As in Example~\ref{ex:prox-descent}, using the Lipschitz continuity of $F$, hence local boundedness of $\partial F$, and using the continuous differentiability of $g$, the sequence of sets $\sum_{i=1}^M \partial F_{i}(\akk{x}\iter\k) (\nabla g(F(x\iter\k)))_i - \partial F_{i}(\akk{x}\iter\k) (\nabla g(F(\akk{x}\iter\k)))_i$ vanishes for $\k\fto{K}\infty$, which implies that the slope $\slp{\fun}(\akk{x}\iter\k)$ vanishes for $k\fto{K}\infty$.
\end{EX}
\begin{EX}[Problem adaptive model function]
Our framework allows for a problem specific adaptation using a combination of Examples~\ref{ex:fbs}, \ref{ex:vm-fbs}, \ref{ex:prox-descent}, and \ref{ex:mm}. Consider the following objective function\footnote{The example is not meant to be meaningful and the model function to be algorithmically the best choice. This example shall demonstrate the flexibility and problem adaptivity of our framework.} $\map{f}{\R^N\times\R^2}{\eR}$ with $a\in \R^N$ and $b\in\R$:
\[
 f(x,z) = f_1(z)+ \ind{[-1,1]^2}(z) + f_2(x) \  \text{with}\ 
 \begin{cases}
    f_1\in \spC2(\R^2)\ \text{is strongly convex}\,; \\
    f_2(x) := \max\{(\scal{a}{x} - b)^2, 1-\exp(-\vnorm{x})\}\,.
 \end{cases}
\]
We define our model function as: $\mfun{(\ak{x},\ak{z})}(x,z) = \bar f_1(z; \ak{z}) +\ind{[-1,1]}(z)+ \bar f_2(x;\ak{x})$ with
\[
 \begin{split}
    \bar f_1(z;\ak{z}) :=&\ f_1(\ak{z}) + \scal{\nabla f_1(\ak{z})}{z-\ak{z}} + \frac12 \scal{z-\ak{z}}{\nabla^2 f_1(\ak{z})(z-\ak{z})} \,;\\
    \bar f_2(x;\ak{x}) :=&\ \max\{(\scal{a_i}{x} - b_i)^2, 1+\exp(-\vnorm{\ak{x}})(\vnorm{x} - \vnorm{\ak{x}}-1)\} \,. \\ 
 \end{split}
\]
The strong convexity of $f_1$ allows for a convex second order approximation with positive definite Hessian. We linearize only the second component of the ``$\max$'' (w.r.t. $\vnorm{x}$) in $f_2$ to obtain a convex approximation that is as close as possible to original function $f_2$. 

As $\nabla f_1$ is Lipschitz continuous on the compact set $[0,1]^2$, the growth function w.r.t. $z$ is of the form $L\vnorm{z-\ak{z}}^2$. Moreover, using $1$-Lipschitz continuity of $\exp$ on $\R_{-}$, we can bound $\abs{\exp(-\vnorm{x}) - \exp(-\vnorm{\ak{x}})(1+\vnorm{\ak{x}} - \vnorm{x})}$ by $\abs{\vnorm{x}-\vnorm{\ak{x}}}^2$ and, using Lipschitz continuity of $\vnorm{x}$, by $\vnorm{x-\ak{x}}^2$. Therefore, the model error is given by a growth function $w(t)=\max\{1,L\}t^2$.
\end{EX}

\subsection{Examples of Bregman functions} \label{subsec:examples-bregman}

Let us explore some of the Bregman functions, that are most important to our applications and show that our assumptions are satisfied.
\begin{EX}[Euclidean Distance] \label{ex:Breg-dist:Euclidean}
The most natural Bregman proximity function is the Euclidean distance 
\[
  \breg[\leg]{x}{\ak{x}}=\frac 12\vnorm{x-\ak{x}}^2\,,
\]
which is generated by the Legendre function $\leg(x) = \frac 12\vnorm{x}^2$. The domain of $\leg$ is the whole space $\X$, which implies that Condition~\ii2 in Theorem~\ref{thm:stationarity-growth} is satisfied for any limit point. Assumption~\ref{ass:Bregman-continuity} is trivial, and for the model functions in Section~\ref{subsec:examples-models}, Assumption~\ref{ass:vanishing-slope} is satisfied, if $x\iter\k - \ak{x}\iter\k\to0$ implies $\nabla \leg(x\iter\k) - \nabla \leg(\ak{x}\iter\k) \to 0$, which is clearly true. Therefore, for the models in Section~\ref{subsec:examples-models} combined with the Euclidean proximity measure, we conclude subsequential convergence to a stationary point.
\end{EX}
\begin{EX}[Variable Euclidean Distance]
A simple but far-reaching extension of Example~\ref{ex:Breg-dist:Euclidean} is the following. Let $\seq[\k\in\N]{A\iter\k}$ be a sequence of symmetric positive definite matrices such that the smallest and largest eigenvalues are in $[c_1,c_2]$ for some $0<c_1<c_2<+\infty$, i.e.  $0< \inf_\k \scal{x}{A\iter\k x} < \sup_\k \scal{x}{A\iter\k x} < +\infty$ for all $x\in \X$. Each matrix $A\iter\k$ induces a metric on $\X$ via the inner product $\scal{x}{A \ak{x}}$ for $x,\ak{x}\in \X$. The induced norm is a Bregman proximity function
\[
\breg[\leg\iter\k]{x}{\ak{x}} = \frac 12\vnorm[A\iter\k]{x-\ak{x}}^2 := \frac 12 \scal{x-\ak{x}}{A\iter\k(x-\ak{x})} \,,
\]
generated analogously to Example~\ref{ex:Breg-dist:Euclidean}. Except the boundedness of the eigenvalues of $\seq[\k\in\N]{A\iter\k}$ there are no other restrictions. All the conditions mentioned in Example~\ref{ex:Breg-dist:Euclidean} are easily verified.

\end{EX}
From now on, we restrict to iteration-independent Bregman distance functions, knowing that we can flexibly adapt the Bregman distance in each iteration. 
\begin{EX}[Boltzmann--Shannon entropy]
The Boltzmann-Shannon entropy is 
\[
  \breg[\leg]{x}{\ak{x}} = \sum_{i=1}^N \big(x\idx{i}(\log(x\idx{i}) - \log(\ak{x}\idx{i})) - (x\idx{i} - \ak{x}\idx{i}) \big)
\]
where $x\idx{i}$ denotes the $i$-th coordinate of $x\in\X$. $\bregmap[\leg]$ is generated by the Legendre function $\leg(x) = \sum_{i=1}^N x\idx{i}\log(x\idx{i})$, which has the domain $[0,+\infty[^N$. Since $\leg$ is additively separable, w.l.o.g., we restrict the discussion to $N=1$ in the following. 

We verify Assumption~\ref{ass:Bregman-continuity}.  Let $\seq[\k\in\N]{x\iter\k}$  and $\seq[\k\in\N]{\ak{x}\iter\k}$ be bounded sequences in $\sint\dom\leg=]0,+\infty[$ with $x\iter\k-\ak{x}\iter\k \to 0$ for $\k\to\infty$. For any convergent subsequence $x\iter\k\to x^*$ as $\k\fto{K}\infty$ for some $K\subset \N$ also $\ak{x}\iter\k \to x^*$ as $\k\fto{K}\infty$ and $x^*\in [0,+\infty[$. Since $\leg$ is continuous on $\scl\dom\leg = [0,+\infty)$ (define $\leg(0) = 0\log(0)=0$), $\breg[\leg]{x\iter\k}{\ak{x}\iter\k}\to 0$ for any convergent subsequence, hence for the full sequence. The same argument shows that the converse implication is also true, hence the Boltzmann-Shannon entropy satisfies Assumption~\ref{ass:Bregman-continuity}.

For the model functions from Section~\ref{subsec:examples-models}, we show that Assumption~\ref{ass:vanishing-slope} holds for $x^*\in \sint\dom\leg$, i.e. $\nabla \leg(x\iter\k) - \nabla \leg(\ak{x}\iter\k)\to 0$ for sequence $\seq[\k\in\N]{x\iter\k}$ and $\seq[\k\in\N]{\ak{x}\iter\k}$ with $x\iter\k\to x^*$ and $x\iter\k - \ak{x}\iter\k\to0$ for $\k\fto{K}\infty$ for some $K\subset \N$. This condition is satisfied, because $\nabla \leg$ is continuous on $\sint\dom\leg$, hence $\lim_{\k\fto{K}\infty} \nabla \leg(x\iter\k) = \lim_{\k\fto{K}\infty} \nabla \leg(\ak{x}\iter\k)=\nabla \leg(x^*)$.

Since $\dom\leg = \scl\dom\leg$, it suffices to verify Condition~\ii3 of Theorem~\ref{thm:stationarity-growth} to guarantee subsequential convergence to a stationary point. For $x^*\in[0,+\infty[$ and a bounded sequence $\seq[\k\in\N]{\akk{y}\iter\k}$ in $\sint\dom\leg$ as in Theorem~\ref{thm:stationarity-growth}, we need to show that $\breg[\leg]{x^*}{\akk{y}\iter\k}\to 0$ as $\k\fto{K} \infty$ for $K\subset \N$ such that $\akk{y}\iter\k\to x^*$ as $\k\fto{K}\infty$. This result is clearly true for $x^*>0$, thanks to the continuity of $\log$. For $x^*=0$, we observe $x^*\log(\akk{y}\iter\k)\to 0$ for $\k\fto{K}\infty$, hence Condition~\ii3 of Theorem~\ref{thm:stationarity-growth} holds, and subsequential convergence to a stationary point is guaranteed.
\end{EX}
\begin{EX}[Burg's entropy]  \label{ex:burg-entropy}
For optimization problems with non-negativity constraint, Burg's entropy is a powerful distance measure. The associated Bregman distance
\[
  \breg[\leg]{x}{\ak{x}} = \sum_{i=1}^N \Bigg( \frac{x\idx{i}}{\ak{x}\idx{i}} - \log\Big(\frac{x\idx{i}}{\ak{x}\idx{i}}\Big) -1\Bigg)
\]
is generated by the Legendre function $\leg(x) = -\sum_{i=1}^N \log(x\idx{i})$ which is defined on the domain $]0,+\infty[^N$. Approaching $0$, the function $\leg$ grows towards $+\infty$. In contrast to the Bregman functions in the examples above, Burg's entropy does not have a Lipschitz continuous gradient, and is therefore interesting for objective functions with the same deficiency. 

W.l.o.g. we consider $N=1$. Assumption~\ref{ass:Bregman-continuity} for two bounded sequences $\seq[\k\in\N]{x\iter\k}$ and $\seq[\k\in\N]{\ak{x}\iter\k}$ in $]0,+\infty[$ reads
\[
  x\iter\k - \ak{x}\iter\k \to 0 \quad \Leftrightarrow\quad    \frac{x\iter\k}{\ak{x}\iter\k} - \log\Big(\frac{x\iter\k}{\ak{x}\iter\k}\Big) \to 1\,,
\]
which is satisfied if the limit points lie in $]0,+\infty[$ since $x\iter\k-\ak{x}\iter\k\to 0 \Leftrightarrow x\iter\k /\ak{x}\iter\k \to 1$ for $\k\fto{K}\infty$ and $\log$ is continuous at $1$. 

For the model functions in Section~\ref{subsec:examples-models}, Assumption~\ref{ass:vanishing-slope} requires $\nabla \leg(x\iter\k) - \nabla \leg(\ak{x}\iter\k)\to 0$ for sequence $\seq[\k\in\N]{x\iter\k}$ and $\seq[\k\in\N]{\ak{x}\iter\k}$ in $\sint\dom\leg$ with $x\iter\k\to x^*$ and $x\iter\k - \ak{x}\iter\k\to0$ for $\k\fto{K}\infty$ for some $K\subset \N$. By continuity, this statement is true for any $x^*>0$. For $x^*=0$, the statement is in general not true. Also Condition~\ii4 in Theorem~\ref{thm:stationarity-growth} can, in general, not be verified. Therefore, if a model functions is complemented with Burg's entropy, then the objective should be continuous on the $\scl\dom\leg$. Stationarity of limit points is obtained, if they lie in $\sint\dom\leg$.
%Condition~\ii4 in Theorem~\ref{thm:stationarity-growth}, which is used to prove $\fun$-attentive convergence, can be verified for Burg's entropy. Let $\seq[\k\in\N]{x\iter\k}$ and $\seq[\k\in\N]{\akk{x}\iter\k}$ be sequences in $\sint\dom\leg=(0,+\infty)$ with $x\iter\k\to x^*=0\in\scl\dom\leg$ and $x\iter\k-\akk{x}\iter\k\to 0$ for $\k\fto{K}\infty$ for some $K\subset \N$. Now, consider any $x\in \sint\dom\leg\cap\dom f$. We make the following estimation:
%\[
%  \breg[\leg]{x}{\akk{x}\iter\k} - \breg[\leg]{x}{x\iter\k} = \frac{x}{\akk{x}\iter\k} - \frac{x}{x\iter\k} - \log\Big(\frac x{\akk{x}\iter\k}\Big) + \log\Big(\frac x{x\iter\k}\Big) = x\Big(\frac 1{\akk{x}\iter\k} - \frac{1}{x\iter\k}\Big) + \log\Big(\frac{\akk{x}\iter\k}{x\iter\k}\Big)
%\]
%The reasoning from before shows that $\log(\akk{x}\iter\k/x\iter\k) \to 0$ as $\k\fto{K}\infty$ and ...\TODO
\end{EX}

% *******************
% >>>>> SECTION <<<<<
% *******************
\section{Applications} \label{sec:numerics}

%Bonettini et al. \cite{BLPP16} have applied Forward--Backward Splitting (FBS), a special setting of our framework (with strongly convex proximity measure), to Student-t regularized image denoising, deblurring under Poisson noise (in a convex setting), and diffusion based image compression. We solve two problems that are beyond their applicability and several dictionary learning formulations. 
We discuss in this section some numerical experiments whose goal is to illustrate the wide applicability of our algorithmic framework. The applicability of our results follows from the considerations in Section~\ref{sec:examples}. Actually, the considered objective functions are all continuous, i.e. Theorem~\ref{thm:stationarity-growth}\ii1 is satisfied.

\subsection{Robust Non-linear Regression} \label{sec:numeric-robust-regression}

We consider a simple non-smooth and non-convex robust regression problem \cite{HRRS86} of the form
\begin{equation} \label{eq:numberic-robust-regression}
  \min_{u:=(a,b)\in \R^{P}\times\R^{P}}\, \sum_{i=1}^M \norm[1]{F_i(u) - y_i} \,, \quad F_i(u) := \sum_{j=1}^P b_j \exp(-a_j x_i) \,,
\end{equation}
where $(x_i,y_i)\in \R\times \R$, $i=1,\ldots,M$ is a sequence of covariate-observation pairs. We assume that $(x_i,y_i)$ are related by $y_i=F_i(u) + n_i$, where $n_i$ is the error term and $u=(a,b)$ are the unknown parameters. We assume that the errors are iid with Laplacian distribution, in which case the data fidelity devised by a maximum likelihood argument is $\ell_1$-norm as used in \eqref{eq:numberic-robust-regression}.

We define model functions by linearizing the inner functions $F_i$ as suggested by the model function in Example~\ref{ex:prox-descent}. Complemented by an Euclidean proximity measure (with $\tau>0$) the convex subproblem \eqref{eq:alg:abstract-relations-A} to be solved inexactly is the following:
\[
  \akk{u} = \argmin_{u\in \R^P\times\R^P}\, \sum_{i=1}^M \norm[1]{\mathcal K_iu - y_i^\diamond} + \frac{1}{2\tau}\vnorm{u - \ak{u}}^2\,,\quad y_i^\diamond := y_i - F(\ak{u}) + \mathcal K_i\ak{u}\,,
\]
where $\map{\mathcal K_i:=DF_i(\ak{u})}{\R^P\times\R^P}{\R}$ is the Jacobian of $F_i$ at the current parameters $\ak{u}$. We solve the (convex) dual problem (cf. \cite{Chambolle04,CDV10}) with warm starting up to absolute step difference $10^{-3}$. 

As mentioned in Remark~\ref{rem:alg-backtracking-subproblems}, backtracking on $\tau$ could be used (cf. ProxDescent \cite{LW16}); denoted \texttt{prox-linear} and \texttt{prox-linear2} in the following. This requires to solve the subproblem for each trial step. This is the bottleneck compared to evaluating the objective. The line search in Algorithm~\ref{alg:line-search} only has to evaluate the objective value. This variant is denoted \texttt{prox-linear-LS} in the following. A representative convergence result in terms of the number of accumulated iterations of the subproblems is shown in Figure~\ref{fig:conv-rob-regr}. For this random example, the maximal noise amplitude is $12.18$, and the maximal absolute deviation of the solution from the ground truth is $0.53$, which is reasonable for this noise level. Algorithm \texttt{prox-linear-LS} requires significantly fewer subproblem iterations than \texttt{prox-linear} and \texttt{prox-linear2}. For \texttt{prox-linear2} the initial $\tau$ is chosen such that initially no backtracking is required.

For large scale problems, frequently solving the subproblems can be prohibitively expensive. Hence, ProxDescent cannot be applied, whereas our algorithm is still practical. 

\tikzset{
  font={\fontsize{10pt}{12}\selectfont}}
\begin{figure}[t]
\begin{center}
\begin{minipage}{0.4\linewidth}
  \begin{tikzpicture}
    \begin{axis}[%
      width=\linewidth,%
      height=5cm,%
      xmin=0,xmax=400000,%
      xlabel=Accumulated subproblem iterations,%
      ylabel=Objective value,%
      legend style={font=\scriptsize}
      ]
    \tiny      
    \addplot[ultra thick,blue,densely dashed] %
        table {figures/prox-linear.dat};
        \addlegendentry{prox-linear};
    \addplot[ultra thick,orange,densely dashed] %
        table {figures/prox-linear_small_tau.dat};
        \addlegendentry{prox-linear2};
    \addplot[ultra thick,green!80!black,densely dotted] %
        table {figures/prox-linear-LS.dat};
        \addlegendentry{prox-linear-LS};
    \end{axis}
  \end{tikzpicture} 
\end{minipage}\hfill
\begin{minipage}{0.55\linewidth}
  \includegraphics[width=0.33\linewidth]{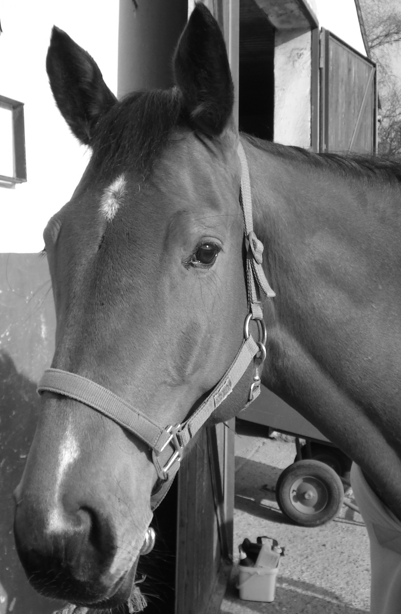}\hspace*{0.3mm}
  \includegraphics[width=0.33\linewidth]{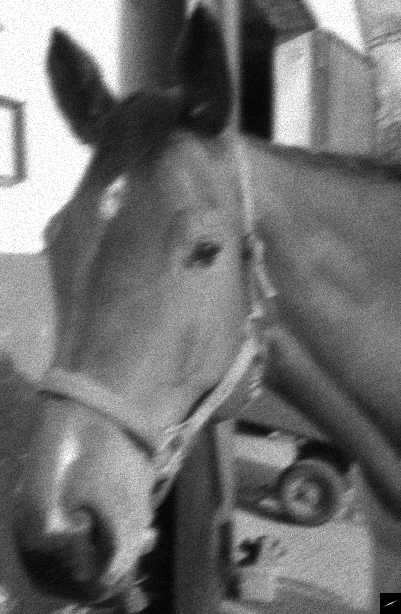}\hspace*{0.3mm}
  \includegraphics[width=0.33\linewidth]{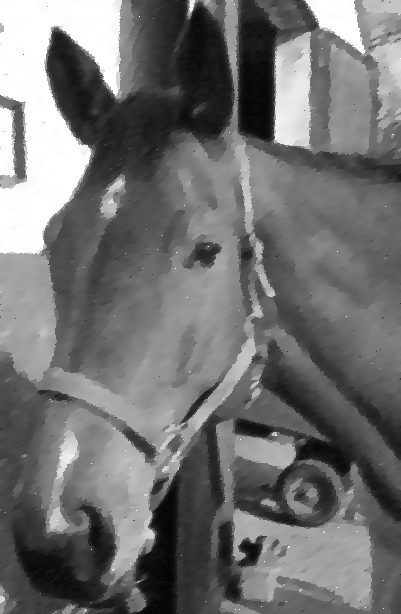}
\end{minipage}
\end{center}
\begin{minipage}[t]{0.4\linewidth}
  \caption{\label{fig:conv-rob-regr}Objective value vs. accumulated number of subproblem iterations for \eqref{eq:numberic-robust-regression}.}
\end{minipage}\hfill
\begin{minipage}[t]{0.55\linewidth}
\caption{\label{fig:posson-img-reconstruction}Deblurring and Poisson noise removal by solving \eqref{eq:poisson-img-prob}. From left to right: clean, noisy, and reconstructed image (PSNR: $25.86$).}
\end{minipage}
\end{figure}

\subsection{Image Deblurring under Poisson Noise}
\newcommand{\Nx}{{n_x}}
\newcommand{\Ny}{{n_y}}
\newcommand{\imgf}{{b}}
\newcommand{\imgu}{{u}}
\newcommand{\imgDiff}{\mathcal D}
\newcommand{\blurOp}{\mathcal A}
\newcommand{\reg}{\varphi}

Let $\imgf\in \R^{\Nx \times \Ny}$ represent a blurry image of size $\Nx\times \Ny$ corrupted by Poisson noise. %, i.e., there exists an image $\imgg$ such that $\imgf = 10^{12} \texttt{imnoise}(\imgg / 10^{12}, \texttt{'Poisson'})$ artificially generated in Matlab. 
Recovering a clean image from $\imgf$ is an ill-posed inverse problem. It is a common problem, for example, in  fluorescence microscopy and optical/infrared astronomy; see \cite{BBDV09} and references therein.  A popular way to solve it is to formulate an optimization problem \cite{ZBZB09} of the form
\begin{equation}\label{eq:poisson-img-prob}
\min_{\imgu\in \R^{\Nx\times \Ny}}\, f(\imgu) := \breg[{KL}]{\imgf}{\blurOp\imgu} + \frac \lambda 2 \sum_{i=1}^\Nx\sum_{j=1}^\Ny \reg(\vnorm{ (\imgDiff \imgu)_{i,j} }^2 ) \,, \quad \st\  \imgu_{i,j} \geq 0\,,
\end{equation}
where $\blurOp$ is a circular convolution (blur) operator. The first term (coined {data term}) in the objective $f$ is the Kullback--Leibler divergence (Bregman distance generated by the Boltzmann--Shannon entropy $x\log(x)$), which, neglecting additive constant terms, is given by 
\[
  f_1(\imgu) := \breg[{KL}]{\imgf}{\blurOp\imgu} := \sum_{i,j} (\blurOp\imgu)_{i,j} - \imgf_{i,j}\log( (\blurOp\imgu)_{i,j})\,,
\]
$f_1$ is well-suited for Poisson noise removal \cite{VSK85}.
The second term (coined {regularization term}) involves a penalty $\map{\reg}{\R^2}{\R}$ applied to spatial finite differences $(\imgDiff \imgu)_{i,j} := ((\imgDiff \imgu)_{i,j}^1 , (\imgDiff \imgu)_{i,j}^2)^\top$ in horizontal direction $(\imgDiff \imgu)_{i,j}^1 := \imgu_{i+1,j} - \imgu_{i,j}$ for all $(i,j)$ with $i<\Nx$, and $0$ otherwise; and vertical direction $(\imgDiff \imgu)_{i,j}^2$ (defined analogously). The function $\reg$ in the regularization is usually chosen to favor \enquote{smooth} images with sharp edges. The relative importance of both the data and regularization terms is weighted by $\lambda>0$. 

For convex penalties $\reg$, algorithms for solving problem \eqref{eq:poisson-img-prob} are available (e.g. primal-dual proximal splitting) provided that $\reg$ is simple (in the sense that its Euclidean proximal mapping can be computed easily). But if one would like to exploit the gradient of $f_1$ explicitly, things become more intricate. The difficulty comes from the lack of global Lipschitz continuity of $\nabla f_1(\imgu)$. A remedy is provided by Bauschke et al. \cite{BBT16}. They have shown that, instead of the global Lipschitz continuity, the key property is the convexity of $L\leg - f_1$ for a Legendre function $\leg$ and sufficiently large $L$, which can be achieved using Burg's entropy $\leg(\imgu) = - \sum_{i,j} \log(\imgu_{i,j})$ (\cite[Lemma~7]{BBT16}).

However, non-convex penalties $\reg$ are known to yield a better solution \cite{GG84,BZ87,MS89}. In this case, the algorithmic framework of Bauschke et al. \cite{BBT16} is not applicable anymore, whereas our framework is applicable. Due to the lack of strong convexity of Burg's entropy also the algorithm of Bonettini et al. \cite{BLPP16} cannot be used. Note that Burg's entropy is strongly convex on bounded subsets of $]0,+\infty[$, however, the subset cannot be determined a priori.

The abstract framework proposed in this paper appears to be the first algorithm with convergence guarantees for solving \eqref{eq:poisson-img-prob} with a smooth non-convex regularizer. 

In our framework, we choose $\varphi: t \in \R^2 \mapsto \log(1+\rho\vnorm{t}^2)$, which is smooth but non-convex. The model functions are defined as in Example~\ref{ex:fbs}. We also use $\leg$ as the Burg's entropy to generate the Bergman proximity function\xspace(see Example~\ref{ex:burg-entropy}). Thus, the subproblems \eqref{eq:alg:abstract-relations-A} which emerge from linearizing the objective $f$ in \eqref{eq:poisson-img-prob} around the current iterate $\ak{\imgu}$
\[
    \akk{\imgu} = \argmin_{\imgu\in \R^{\Nx\times \Ny}} \scal{\imgu-\ak{\imgu}}{\nabla f(\ak{\imgu})} + \frac{1}{\tau} \sum_{i,j}\left( \frac{\imgu_{i,j}}{\ak{\imgu}_{i,j}} -  \log\Big(\frac{\imgu_{i,j}}{\ak{\imgu}_{i,j}}\Big) \right) 
\]
can be solved exactly in closed-form $\akk{\imgu}_{i,j} = \ak{\imgu}_{i,j}/(1+\tau (\nabla f(\ak{\imgu}))_{i,j} \ak{\imgu}_{i,j})$ for all $i,j$. A result for the successful Poisson noise removal and deblurring is shown in Figure~\ref{fig:posson-img-reconstruction}.

\subsection{Structured Matrix Factorization} \label{sec:numeric:dictionary-learning}

Structured matrix factorization problems are crucial in data analysis. It has many applications in various areas including blind deconvolution in signal processing, clustering, source separation, dictionary learning, etc.. There is a large body of literature on the subject and we refer to e.g.~\cite{Cichocki09,Chaudhuri14,Starck2015sparse,XuDL17} and references therein for a comprehensive account.

\paragraph{The problem.}
Given a data matrix $A\in \R^{M\times N}$ whose $N$ $M$-dimensional columns are the data vectors. The goal is to find two matrices $U\in \R^{M\times K}$ and $Z\in \R^{K\times N}$ such that 
\[
  A = U Z + Q\,,
\]
where $Q\in \R^{M\times N}$ accounts for an unknown error. The matrices $U$ and $Z$ (called also factors) enjoy features arising in a specific application at hand (see more below).\\
%For descriptive dictionary a representation matrix 

%In fact, this problem can be seen as a matrix factorization problem which goes beyond dictionary learning. Indeed, most of the subsequent ideas can be applied similarly to, e.g., blind deconvolution or non-negative matrix factorization. Thus, by no means, we intend to improve the state-of-the-art in modelling of the dictionary learning problem, which is outside the scope of the present paper. Rather, we show how the versatility of our abstract algorithmic framework allows to easily adapt to various constraints on $U$ and $Z$, which are just made concrete for the dictionary learning problem (the interested reader may refer to e.g.~\cite{XuDL17} for a recent survey on dictionary learning for face recognition).

%\paragraph{Generic optimization problem.} %Depending on the goal, we are interested in dictionaries with different properties. We formulate different optimization problems, all of the following generic form:
%The following generic optimization problems allows us to recover dictionaries with different properties:
To solve the matrix factorization problem, we adopt the optimization approach and we consider the
non-convex and non-smooth minimization problem
\begin{equation} \label{eq:numeric:dictionary-learning-model}
  \min_{U\in \mathcal U, Z\in \mathcal Z}\, f(A,UZ) + \lambda g(Z) \,,\quad f(A,UZ) := \frac 12\norm[F]{A - UZ}^2 \,.
\end{equation}
The term $f(A,UZ)$ stands for proximity function that measures fidelity of the approximation of $A$ by the product $UZ$ of the two factors. We here focus on the classical case where the fidelity is measured via the Frobenius norm $\norm[F]\cdot$, but other data fidelity measures can also be used just as well in our framework, such as divergences (see \cite{Cichocki09} and references therein). The sets $\mathcal U$, $\mathcal Z$, which are non-empty closed and convex, and the function $g \in \Gamma_0$ are used to capture specific features of the matrices $U$ and $Z$ arising in a specific application as we will exemplify shortly. The influence of $g$ is weighted by the parameter $\lambda > 0$. 

Many (if not most) algorithms to solve the matrix factorization problem~\eqref{eq:numeric:dictionary-learning-model} are based on Gauss-Seidel alternating minimization with limited convergence guarantees~\cite{Cichocki09,Chaudhuri14,Starck2015sparse}\footnote{For very specific instances, a recent line of research proposes to lift the problem to the space of low-rank matrices, and then use convex relaxation and computationally intensive conic programming that are only applicable to small-dimensional problems; see, e.g.,~\cite{Romberg14} for blind deconvolution.}. The PALM algorithm proposed recently by Bolte et al.~\cite{BST14}, was designed specifically for the structure of the optimization problem~\eqref{eq:numeric:dictionary-learning-model}. It can then be successfully applied to solve instances of such a problem with provably guaranteed convergence under some assumptions including the Kurdyka-{\L}ojasiewicz property. However, though it can handle non-convex constraint sets and functions $g$, it does not allow to incorporate Bregman proximity functions. 

In the following, we show how our algorithmic framework can be applied to a broad class of matrix factorization instances. In particular, a distinctive feature of our algorithm is that it can naturally and readily accommodate for different Bregman proximity functions and it has no restrictions on the choice of the step size parameters (except positivity). A descent is enforced in the line search step, which follows the proximal step.

\paragraph{A generic algorithm.} We apply Algorithm~\ref{alg:abstract-alg} to solve this problem, where the model functions are chosen to linearize the data fidelity function $f(A,UZ)$, according to Example~\ref{ex:fbs}. The convex subproblems to be solved in the algorithm have the following form:
\[
  \begin{split}
  (\akk{U},\akk{Z}) = \argmin_{U\in \mathcal U, Z\in \mathcal Z}\, \lambda g(Z) +&\ \scal{Z-\ak{Z}}{\ak{U}^\top(\ak{U}\ak{Z}-A)}_F + \breg[\leg_Z]{Z}{\ak{Z}}  \\
                                                                         +&\ \scal{U - \ak{U}}{(\ak{U}\ak{Z}-A)\ak{Z}^\top}_F + \breg[\leg_U]{U}{\ak{U}} \,
  \end{split}
\]
where $\scal{\cdot}{\cdot}_F$ stands for the Frobenius inner product. The Bregman proximity functions $\breg[\leg_Z]{\cdot}{\cdot}$ and $\breg[\leg_U]{\cdot}{\cdot}$ provide the flexibility to handle a variety of constraint sets $\mathcal U$ and $\mathcal Z$. In the following, we list different choices for the constraint sets and explain how to incorporate them into the optimization procedure. Due to the structure of the optimization problem, the variables $U$ and $Z$ can be handled separately. The only coupling is the data fidelity function $f$, which is linearized and therefore easy to incorporate.

\paragraph{Examples of constraints $\mathcal U$.} There are many possible choices for the set $\mathcal U$ depending on the application at hand.
\begin{itemize}
\item \emph{Unconstrained case}:
  \[
    \mathcal U_1 = \R^{M\times K} \,.
  \]
  In the unconstrained case, a suitable Bregman proximity function is given by the Euclidean distance $\breg[\leg_U]{U}{\ak{U}} = \frac{1}{2\tau_U}\norm[F]{U-\ak{U}}^2$ with step size parameter $\tau_U$. The resulting update step with respect to the dictionary $U$ is a gradient descent step.
  
\item \emph{Zero-mean and normalization}:
\[
  \mathcal U_2 = \set{U\in \R^{M\times K}\setsep 	\forall j \colon \sum_{i=1}^M U_{i,j}^2 \leq 1 \,, 
                                                 	\forall j\geq 2\colon \sum_{i=1}^M U_{i,j} = 0 }\,.
\]
%\JF{The normalization on the unit sphere makes $\mathcal U_2$, hence the model function, non-convex. Moreover, the projector is not single valued everywhere. I changed the sphere to the ball.} 

This choice of the constraint set leads to a natural normalization of the columns of $U$ that removes the scale ambiguity due to bilinearity. This choice is very classical in dictionary learning, see, e.g.,~\cite{XuDL17}. As in dictionary learning, the average of the first column may not be enforced to be zero, in order to allow the first column to absorb the mean value of the data points. 

By separability of $\mathcal U_2$, the Euclidean projection onto it is simple. This projector is column-wise achieved by subtracting the mean, and then projecting the result onto the Euclidean unit ball. Thus we advocate $\breg[\leg_U]{U}{\ak{U}} = \frac1{2\tau_U}\norm[F]{U-\ak{U}}^2$ with step size parameter $\tau_U$. In turn, the subproblem with respect to $U$ amounts to a projected gradient descent step. 

\item \emph{Non-negativity and normalization}:
\[
  \mathcal U_3 = \set{U\in \R^{M\times K}\setsep \forall j \colon \sum_{i=1}^M U_{i,j} = 1\,,\ 
                                                 \forall i,j\colon U_{i,j} \geq 0 } \,.
\]
This choice is adopted in non-negative matrix factorization (NMF)~\cite{Lee99}. 
%Using a dictionary with non-negative entries, the \enquote{properties} of the data points are clearly separated. A negative representation coefficient shows that the associated property (the dictionary atom) is a negative characterization of this data point, and a positive representation coefficient shows that the data point contains the property represented by this coefficient. In order to describe all data points with positive properties only, the non-negativity constraint of the dictionary can be complemented with non-negativity of the representation coefficients. 
The constraint set $\mathcal U_3$ is column-wise a unit simplex constraint. This constraint can be conveniently  handled by choosing $\breg[\leg_U]{U}{\ak{U}} = \frac{1}{\tau_U}\sum_{i,j} U_{i,j}(\log(U_{i,j}) - \log(\ak{U}_{i,j})) - U_{i,j} + \ak{U}_{i,j}$, which is the Bregman function generated by the entropy $\leg_U(U) = \frac{1}{\tau_U} \sum_{i,j} U_{i,j}\log(U_{i,j})$ with step size parameter $\tau_U$. This is a more natural choice than the Euclidean proximity distance. Indeed, the update step with respect to $U$ results in 
\[
\akk{U}_{i,j} = \frac{\ak{U}_{i,j}\exp(-\tau_U (C_U)_{i,j})}{\sum_{p=1}^{M} \ak{U}_{p,j}\exp(-\tau_U (C_U)_{p,j})} \quad \forall i=1,\ldots, M\,;\ \forall j=1,\ldots,K\,,
\]
where we use the shorthand notation $C_U:=\nabla_U f(A,\ak{U}\ak{Z})=\ak{U}^\top(\ak{U}\ak{Z}-A)$ for the partial gradient of $f$ with respect to $U$. The exponential function is applied entry-wise, hence naturally preserving positivity. Note that the Euclidean projector onto $\mathcal U_3$ necessitates to compute the projector on the simplex which can be achieved with sorting~\cite{Michelot86}.
\end{itemize}

\paragraph{Examples of constraints $\mathcal Z$.}
There are also several possible choices for the set $\mathcal Z$ and regularizing function $g$ depending on the application at hand.
\begin{itemize}
  \item \emph{Unconstrained case}: 
  \[
    \mathcal Z_1 = \R^{K\times N} \quad\text{and}\quad g(Z) = 0\,.
  \]
  This case can be handled using a gradient descent step, analogously to the related update step with the constraint set $\mathcal U_1$.
  \item \emph{Non-negativity}:
  \[
    \mathcal Z_2 = \set{Z\in \R^{K\times N}\setsep \forall i,j\colon Z_{i,j} \geq 0 } \quad\text{and}\quad g(Z) = 0\,.
  \]
  This constraint is used in conjunction with $\mathcal U_3$ in NMF. It can be handled either with a Euclidean proximity function (which amounts to projecting on the non-negative orthant), or via a Bregman proximity function $\breg[\leg_Z]{Z}{\ak{Z}}$ generated by the Boltzmann--Shannon entropy ($\leg_Z(Z) = \frac{1}{\tau_Z}\sum_{i,j}Z_{i,j}\log(Z_{i,j})$) or, alternatively, Burg's entropy \linebreak ($\leg_Z(Z) = -\frac{1}{\tau_Z}\sum_{i,j}\log(Z_{i,j})$), with step size parameter $\tau_Z$. The update with respect to $Z$ then reads
\[
\akk{Z}_{i,j} = \ak{Z}_{i,j}\exp(-\tau_Z (C_Z)_{i,j}) \quad \forall i=1,\ldots, K\,;\ \forall j=1,\ldots,N\,,
\]
where we use the shorthand notation $C_Z:=\nabla_Z f(A,\ak{U}\ak{Z})=(\ak{U}\ak{Z}-A)\ak{Z}^\top$ for the partial gradient of $f$ with respect to $Z$.
 
  \item \emph{Sparsity constraints}: 
  \[
    \mathcal Z_3 = \R^{K\times N} \quad \text{and}\quad
      g(Z) =  \norm[1]{Z} \,.
  \]

The introduction of sparsity has been of prominent importance in several matrix factorization problems, including dictionary learning~\cite{OF97}, NMF~\cite{Hoyer04} \footnote{Strictly speaking, $\mathcal Z_3$ should be the non-negative orthant for sparse NMF. But this does not change anything to our discussion since computing the Euclidean proximal mapping of the $\ell_1$ norm restricted to the non-negative orthant is easy.} and source separation~\cite{Starck2015sparse}. The Euclidean proximal mapping of the $\ell_1$-norm is the entry-wise soft-thresholding, hence giving the update step with respect to $Z$ as
\[
\akk{Z}_{i,j} = \max\{0,1-\lambda\tau_Z/|\ak{Z}_{i,j}-\tau_Z (C_Z)_{i,j}|\}(\ak{Z}_{i,j}-\tau_Z (C_Z)_{i,j}) \quad \forall i=1,\ldots, K\,;\ \forall j=1,\ldots,N\, .
\]

   \item \emph{Low rank constraint}: 
   \[
    \mathcal Z_3 = \R^{K\times N} \quad \text{and}\quad
      g(Z) =  \norm[*]{Z} \,.
  \]
   
The nuclear norm or 1-Schatten norm $\norm[*]{Z}$ is the sum of the singular values. It is known to be the tightest convex relaxation to the rank and was shown to promote low rank solutions~\cite{recht2010guaranteed}. Such a regularization would be useful in the situation where columns of $A$ are (to a good approximation) clustered on a few linear subspaces spanned by the columns of $U$, i.e. the columns of $A$ can be explained by columns of $U$ from the same subspace (\enquote{cluster}).
   
The Euclidian proximal mapping of the nuclear norm is the soft-thresholding applied to the singular values. In turn, the update step with respect to $Z$ reads
\[
\akk{Z}_{i,j} = W\diag((\max\{0,1-\lambda\tau_Z/\sigma_i\}\sigma_i)_i)V^\top ,
\]
where $W$, $V$ are respectively the matrices of left and right singular vectors of $\ak{Z}-\tau_Z C_Z$, and $\sigma$ is the associated vector of singular values.
\end{itemize}

\section{Conclusions}

We have presented an algorithmic framework, that unifies the analysis of several first order optimization algorithms in non-smooth non-convex optimization such as Gradient Descent, Forward--Backward Splitting, ProxDescent, and many more. The algorithm combines sequential Bregman proximal minimization of model functions, which is the key concept for the unification, with an Armijo-like line search strategy. The framework reduces the difference between algorithms to the model approximation error measured by a growth function. For the developed abstract algorithmic framework, we establish subsequential convergence to a stationary point and demonstrate its flexible applicability in several difficult inverse problems from machine learning, signal and image processing.

% ************************
% >>>>> bibliography <<<<<
% ************************
{\small
\bibliographystyle{ieee}
\bibliography{ochs}

\begin{thebibliography}{10}\itemsep=-1pt

\bibitem{Romberg14}
A.~Ahmed, B.~Recht, and J.~Romberg.
\newblock Blind deconvolution using convex programming.
\newblock {\em IEEE Transactions on Information Theory}, 60(3):1711--1732,
  2014.

\bibitem{ABS13}
H.~Attouch, J.~Bolte, and B.~Svaiter.
\newblock Convergence of descent methods for semi-algebraic and tame problems:
  proximal algorithms, forward--backward splitting, and regularized
  {G}auss--{S}eidel methods.
\newblock {\em Mathematical Programming}, 137(1-2):91--129, 2013.

\bibitem{BB97}
H.~Bauschke and J.~Borwein.
\newblock Legendre functions and the method of random {Bregman} projections.
\newblock {\em Journal of Convex Analysis}, 4(1):27--67, 1997.

\bibitem{BBC01}
H.~Bauschke, J.~Borwein, and P.~Combettes.
\newblock Essential smoothness, essential strict convexity, and {L}egendre
  functions in {B}anach spaces.
\newblock {\em Communications in Contemporary Mathematics}, 3(4):615--647, Nov.
  2001.

\bibitem{BBC03}
H.~Bauschke, J.~Borwein, and P.~Combettes.
\newblock Bregman monotone optimization algorithms.
\newblock {\em SIAM Journal on Control and Optimization}, 42(2):596--636, Jan.
  2003.

\bibitem{BBT16}
H.~H. Bauschke, J.~Bolte, and M.~Teboulle.
\newblock A descent lemma beyond {L}ipschitz gradient continuity: First-order
  methods revisited and applications.
\newblock {\em Mathematics of Operations Research}, 42(2):330--348, Nov. 2016.

\bibitem{BC11}
H.~H. Bauschke and P.~L. Combettes.
\newblock {\em Convex analysis and monotone operator theory in Hilbert spaces}.
\newblock Springer, 2011.

\bibitem{BBDV09}
M.~Bertero, P.~Boccacci, G.~Desider\`{a}, and G.~Vicidomini.
\newblock Image deblurring with {Poisson} data: from cells to galaxies.
\newblock {\em Inverse Problems}, 25(12):123006, 2009.

\bibitem{BZ87}
A.~Blake and A.~Zisserman.
\newblock {\em Visual Reconstruction}.
\newblock {MIT} Press, Cambridge, {MA}, 1987.

\bibitem{BST14}
J.~Bolte, S.~Sabach, and M.~Teboulle.
\newblock Proximal alternating linearized minimization for nonconvex and
  nonsmooth problems.
\newblock {\em Mathematical Programming}, 146(1-2):459--494, 2014.

\bibitem{BLPP16}
S.~Bonettini, I.~Loris, F.~Porta, and M.~Prato.
\newblock Variable metric inexact line-search based methods for nonsmooth
  optimization.
\newblock {\em SIAM Journal on Optimization}, 26(2):891--921, Jan. 2016.

\bibitem{Bregman67}
L.~M. Bregman.
\newblock The relaxation method of finding the common point of convex sets and
  its application to the solution of problems in convex programming.
\newblock {\em {USSR} Computational Mathematics and Mathematical Physics},
  7(3):200--217, 1967.

\bibitem{Burg72}
J.~Burg.
\newblock The relationship between maximum entropy spectra and maximum
  likelihood spectra.
\newblock {\em Geophysics}, 37(2):375--376, Apr. 1972.

\bibitem{Chambolle04}
A.~Chambolle.
\newblock An algorithm for total variation minimization and applications.
\newblock {\em Journal of Mathematical Imaging and Vision}, 20:89--97, 2004.

\bibitem{Chaudhuri14}
S.~Chaudhuri, R.~Velmurugan, and R.~Rameshan.
\newblock {\em Blind Image Deconvolution}.
\newblock Springer, 2014.

\bibitem{CT93}
G.~Chen and M.~Teboulle.
\newblock Convergence analysis of proximal-like minimization algorithm using
  bregman functions.
\newblock {\em SIAM Journal on Optimization}, 3:538--543, 1993.

\bibitem{Cichocki09}
A.~Cichocki, R.~Zdunek, A.~Phan, and S.~Amari.
\newblock {\em Nonnegative Matrix and Tensor Factorizations: Applications to
  Exploratory Multi-Way Data Analysis and Blind Source Separation}.
\newblock Wiley,, New York, 2009.

\bibitem{CDV10}
P.~Combettes, D.~D{\~{u}}ng, and B.~V{\~{u}}.
\newblock Dualization of signal recovery problems.
\newblock {\em Set-Valued and Variational Analysis}, 18(3-4):373--404, Dec.
  2010.

\bibitem{DIL16}
D.~Drusvyatskiy, A.~D. Ioffe, and A.~S. Lewis.
\newblock Nonsmooth optimization using {Taylor}-like models: error bounds,
  convergence, and termination criteria.
\newblock {\em ArXiv e-prints}, Oct. 2016.
\newblock arXiv: 1610.03446.

\bibitem{DL16}
D.~Drusvyatskiy and A.~S. Lewis.
\newblock Error bounds, quadratic growth, and linear convergence of proximal
  methods.
\newblock {\em ArXiv e-prints}, Feb. 2016.
\newblock arXiv:1602.06661.

\bibitem{GG84}
S.~Geman and D.~Geman.
\newblock Stochastic relaxation, {G}ibbs distributions, and the {B}ayesian
  restoration of images.
\newblock {\em IEEE Transactions on Pattern Analysis and Machine Intelligence},
  6:721--741, 1984.

\bibitem{HRRS86}
F.~R. Hampel, E.~M. Ronchetti, P.~J. Rousseeuw, and W.~A. Stahel.
\newblock {\em Robust Statistics: The Approach Based on Influence Functions}.
\newblock {MIT} Press, Cambridge, {MA}, 1986.

\bibitem{Hoyer04}
P.~Hoyer.
\newblock Non-negative matrix factorization with sparseness constraints.
\newblock {\em J. Mach. Learn. Res.}, 5:1457--1469, 2004.

\bibitem{Lee99}
D.~Lee and H.~Seung.
\newblock Learning the part of objects from nonnegative matrix factorization.
\newblock {\em Nature}, 401:788--791, 1999.

\bibitem{LW16}
A.~Lewis and S.~Wright.
\newblock A proximal method for composite minimization.
\newblock {\em Mathematical Programming}, 158(1-2):501--546, July 2016.

\bibitem{LM79}
P.~L. Lions and B.~Mercier.
\newblock Splitting algorithms for the sum of two nonlinear operators.
\newblock {\em SIAM Journal on Applied Mathematics}, 16(6):964--979, 1979.

\bibitem{Marquardt63}
D.~Marquardt.
\newblock An algorithm for least-squares estimation of nonlinear parameters.
\newblock {\em Society for Industrial and Applied Mathematics}, 11:431--441,
  1963.

\bibitem{Michelot86}
C.~Michelot.
\newblock A finite algorithm for finding the projection of a point onto the
  canonical simplex of $\mathbb{R}^n$.
\newblock {\em J. Optim. Theory Appl.}, 50:195--200, 1986.

\bibitem{MS89}
D.~Mumford and J.~Shah.
\newblock Optimal approximations by piecewise smooth functions and associated
  variational problems.
\newblock {\em Communications on Pure and Applied Mathematics}, 42:577--685,
  1989.

\bibitem{Nguyen17}
Q.~Nguyen.
\newblock Forward--{Backward} {Splitting} with {Bregman} {Distances}.
\newblock {\em Vietnam Journal of Mathematics}, pages 1--21, Jan. 2017.

\bibitem{Noll13}
D.~Noll.
\newblock Convergence of non-smooth descent methods using the
  {Kurdyka}--{{\L}ojasiewicz} inequality.
\newblock {\em Journal of Optimization Theory and Applications},
  160(2):553--572, Sept. 2013.

\bibitem{NPA08}
D.~Noll, O.~Prot, and P.~Apkarian.
\newblock A proximity control algorithm to minimize nonsmooth and nonconvex
  functions.
\newblock {\em Pacific Journal of Optimization}, 4(3):571--604, 2008.

\bibitem{ODBP15}
P.~Ochs, A.~Dosovitskiy, T.~Brox, and T.~Pock.
\newblock On iteratively reweighted algorithms for nonsmooth nonconvex
  optimization in computer vision.
\newblock {\em SIAM Journal on Imaging Sciences}, 8(1):331--372, 2015.

\bibitem{OF97}
B.~Olshausen and D.~Field.
\newblock Sparse coding with an overcomplete basis set: A strategy employed by
  {V1}?
\newblock {\em Vision Research.}, 37, 1996.
\newblock 3311--3325.

\bibitem{recht2010guaranteed}
B.~Recht, M.~Fazel, and P.~A. Parrilo.
\newblock Guaranteed minimum-rank solutions of linear matrix equations via
  nuclear norm minimization.
\newblock {\em SIAM review}, 52(3):471--501, 2010.

\bibitem{Rock98}
R.~T. Rockafellar and R.-B. Wets.
\newblock {\em Variational Analysis}, volume 317.
\newblock Springer Berlin Heidelberg, Heidelberg, 1998.

\bibitem{Starck2015sparse}
J.-L. Starck, F.~Murtagh, and J.~Fadili.
\newblock {\em Sparse image and signal processing: wavelets, curvelets,
  morphological diversity}.
\newblock Cambridge University Press, 2nd edition, 2015.

\bibitem{VSK85}
Y.~Vardi, L.~Shepp, and L.~Kaufman.
\newblock A statistical model for positron emission tomography.
\newblock {\em Journal of the American Statistical Association}, 80(389):8--20,
  1985.

\bibitem{XuDL17}
Y.~Xu, Z.~Li, J.~Yang, and D.~Zhang.
\newblock A survey of dictionary learning algorithms for face recognition.
\newblock {\em IEEE Access}, 5:8502--8514, 2017.

\bibitem{ZBZB09}
R.~Zanella, P.~Boccacci, L.~Zanni, and M.~Bertero.
\newblock Efficient gradient projection methods for edge-preserving removal of
  {Poisson} noise.
\newblock {\em Inverse Problems}, 25(4), 2009.

\end{thebibliography}
}

\end{document}